\documentclass[a4paper,leqno]{article}

\usepackage{amsmath}
\usepackage{amsfonts}
\usepackage{enumerate}
\usepackage{theorem}

\theoremstyle{plain}
\newtheorem{teo}{Theorem}[section]
\newtheorem{lem}[teo]{Lemma}

\newtheorem{prop}[teo]{Proposition}
\newtheorem{defin}[teo]{Definition}

{\theorembodyfont{\rmfamily}
\newtheorem{oss}[teo]{Remark}
}

\renewcommand{\eqref}[1]{\textnormal{(\ref{#1})}}

\numberwithin{equation}{section}

\newcommand{\cvd}{\hfill$\square$}

\title{Stable Determination of the Discontinuous Conductivity Coefficient of a Parabolic Equation
\footnotemark[1]
}

\author{
Michele Di Cristo\footnotemark[2],\hspace{1em}Sergio Vessella\footnotemark[3]\\
}

\date{}

\begin{document}

\maketitle
\footnotetext[1]{Work supported by MIUR, PRIN n.~2006014115.}
\footnotetext[2]{Dipartimento di Matematica,
Politecnico di Milano, Italy.\\ E-mail: \texttt{michele.dicristo@polimi.it}}
\footnotetext[3]{Dipartimento di Matematica per le Decisioni,
Universit\`a degli Studi di Firenze, Italy. E-mail: \texttt{sergio.vessella@dmd.unifi.it}}

\setcounter{section}{0}
\setcounter{secnumdepth}{2}

\begin{abstract}
We deal with the problem of determining a time varying inclusion within
a thermal conductor. In particular we study the continuous
dependance of the inclusion from the
Dirichlet--to--Neumann map.
Under a priori regularity assumptions on the
unknown defect we establish logarithmic stability estimates.

\noindent\textbf{AMS 2000 Mathematics Subject Classification}
Primary 35R30. Secondary 35B60, 33C90.

\noindent \textbf{Keywords} inverse problems, stability, parabolic
equations, unique continuation.
\end{abstract}

\section{Introduction}
In this paper we study the stability issue for the inverse problem of
recovery the discontinuous conductivity coefficient of a parabolic equation
from infinitely many boundary mesurements.

First let us give a coarse formulation of the problem which we are going to
study. Let $T$ be a given positive number. Let $\Omega $ be a bounded domain
of $\mathbb{R}^{n}$, $n\geq 2$, with a sufficiently smooth boundary and let $Q$ be a
domain contained in $\Omega \times \left( 0,T\right) $. Assume that for
every $\tau \in \left( 0,T\right) $ the intersection $D(\tau)$ of $Q$ with the
hyperplane $t=\tau $ is a nonempty set and $\Omega\setminus \overline{D(\tau)}$ is connected
and denote by $k$, $k\neq 1$ a positive
constant. Let $u$ be the weak solution to the following parabolic
initial-boundary value problem
$$\left\{\begin{array}{ll}
\partial_tu-\mathrm{div}((1+(k-1)\chi_Q)\nabla u)=0 &\text{in }\Omega \times(0,T),\\[1mm]
u(\cdot,0)=0 &\text{in }\overline{\Omega},\\[1mm]
u=g &\text{on }\partial\Omega\times (0,T),
\end{array}\right.$$
where $g$ is a prescribed function on $\partial \Omega \times \left(
0,T\right) $. The inverse problem we are addressing to is to determine the
region $Q$ when infinitely many boundary measurements $\left\{ g,\dfrac{%
\partial u}{\partial \nu }_{\mid \partial \Omega \times \left( 0,T\right)
}\right\} $ are available. The problem formulated above arises in
nondestructive testing evaluation (\cite{Ca-Mo}, \cite{Pa-La-Al}).

A uniqueness result for the problem introduced above has been proven in 1997
by Elayyan and Isakov \cite{El-Is}. The main tools on which the approach of
[El-Is] is based are the Runge approximation property and the use of
solutions with Green's function type singularity. For the nonconstructive
character of the Runge property, such an argument does not seem suitable for
our purpose of proving an accurate stability estimate of $Q$ under some a
priori information. Thus, along the line of previous elliptic and parabolic
inverse problems \cite{Al-DC}, \cite{Al-Ve}, \cite{Ve}, we abandon such an approach and we
choose to use arguments based on quantitative estimates of unique
continuation \cite{Al-Be-Ro-Ve}, \cite{DC-Ro-Ve}, \cite{Mo-Ro, Mo-Ro2}.
We also make use of singular solutions of Green's type, but
more quantitative information are necessary in order to obtain stability
estimates. In particular we need an accurate study of the asymptotic
behaviour when the singularity gets close to the interface $\partial Q$.

In the present paper we prove that, under mild a priori assumption on the
topology and the regularity of $Q$, such an inclusion depends continuously
on the boundary data with a rate of continuity of logarithmic type (see
Theorem \ref{maintheorem} for a the precise statement of the result).
In the context of elliptic inverse problems, it has been shown that
logarithmic stability estimates are optimal (\cite{DC-Ro}).
For parabolic inverse problems with unknown boundaries (and the whole Dirichlet--to--Neumann map)
examples showing that the continuous dependance can be at most of logarithmic type,
have been obtained in \cite{DC-Ro-Ve} and \cite{Ve}.
Their proofs work in our situation as well. Namely in such papers the
limit cases $k=+\infty$ and $k=0$ are considered. Everything remains basically
the same in the intermediate situation with $0<k<\infty$, $k\neq1$.

A crucial tool to obtain the logarithmic stability estimate
is connected with a precise evaluation of smallness propagation based on the
two-sphere one-cylinder inequality for solution to parabolic equations
\cite{Es-Fe-Ve}, \cite{Ve} (Theorem 3.10, in the
present paper). Indeed, roughly speaking, such an inequality allows us to
approach the boundary of the inclusion in any slice of time.

Finally we wish to mention here papers of Daido, Kang and Nakamura \cite{Da-Ka-Na}
and, more recently, Isakov, Kim and Nakamura \cite{IKN}
which are strictly related to the present one.
In \cite{IKN}, the authors consider a similar inverse parabolic problem
of detecting an inclusion, that does not depend on the time, by mean of infinitely
boundary measurements and provide a reconstruction procedure to identify it.

In this paper we have decided not to deal with the case $n=1$. Let us just
observe that such a case is easier and it can be treated essentially with similar
assumptions regarding the topology of the set $\Omega\setminus\overline{D(t)}$
(see Remark \ref{n=1} below).

The plan of the paper is the following.
In {\em Section 2} we state our main result. We first give the notations and definitions
we need throughout the paper ({\em Subsection 2.1}) and then in {\em Subsection 2.2} we
state the hypothesis and the stability theorem (Theorem \ref{maintheorem}).
In {\em Section 3} we provide a proof of Theorem \ref{maintheorem}.
We derive first some identity which will be the starting point of our proof.
Then we give some auxiliary result concerning the Hausdorff distance
(Proposition \ref{dh}), \ref{dm2} and \ref{grafrel}),
fundamental solutions (Proposition \ref{5-16Tpr})
and unique continuation properties (Theorem \ref{2sfe1cil}).
Afterward, using the assumptions on the regularity of the inclusion,
we derive some further property related to the distance of two inclusions
and state Proposition \ref{prefb} which provides lower bounds for the solution of the problem.
Finally we prove Theorem \ref{maintheorem}.
Proof of auxiliary propositions are given in {\em Section 4}.
Proposition \ref{dm2} is proven in {\em Subsection 4.1}.
In {\em Subsection 4.2} we prove Proposition \ref{5-16Tpr}
and we also give an asymptotic estimate for the fundamental solution (Theorem \ref{STIMA ASINTOTICA})
which will be used in the next {\em Subsection 4.3} for the proof of Proposition \ref{prefb}.

\section{The Main Result}
\label{mr}

\subsection{Notations and Definitions}\label{nd}
For every $x\in\mathbb{R}^n$, with $n\geq2$, $x=(x_1,\dots,x_n)$,
we shall set $x=(x',x_n)$, where $x'=(x_1,\dots,x_{n-1})\in\mathbb{R}^{n-1}$.
We shall use $X=(x,t)$ to denote a point in $\mathbb{R}^{n+1}$,
where $x\in\mathbb{R}^{n}$ and $t\in\mathbb{R}$.
For every $x\in\mathbb{R}^n$ and $X=(x,t)\in\mathbb{R}^{n+1}$,
we shall set
$$|x|=\left(\sum_{i=1}^nx_i^2\right)^{1/2},\qquad|X|=\left(|x|^2+|t|\right)^{1/2}.$$
Let $r$ be a positive number. For $x_0\in\mathbb{R}^n$ we shall denote
$B_r(x_0)=\{x\in\mathbb{R}^n\,:\,|x-x_0|<r\}$ and
$B'_r(x'_0)=\{x'\in\mathbb{R}^{n-1}\,:\,|x'-x'_0|<r\}$.
We generally set $B_r=B_r(0)$ and $B'_r=B'_r(0)$.
We denote by $B_r^+(0)=\{x\in B_r\,:\, x_n>0\}$
and $B_r^-(0)=\{x\in B_r\,:\, x_n<0\}$.
For a point $X_0=(x_0,t_0)\in\mathbb{R}^{n+1}$ we shall set
$Q_r(X_0)= B_r(x_0)\times(t_0-r^2,t_0)$.

Given a subset $A$ of $\mathbb R^n$, we shall denote by
\begin{eqnarray*}
&&[A]_{\varepsilon}=\left\{x\in\mathbb R^n\,:\,\mathrm{dist}(x,A)\leq\varepsilon\right\},\\[2mm]
&&(A)_{\varepsilon}=\left\{x\in A\,:\,\mathrm{dist}(x,\partial A)>\varepsilon\right\},\\[2mm]
&&[\partial A]_{\varepsilon}=\left\{x\in\mathbb R^n\,:\,\mathrm{dist}(x,\partial A)\leq\varepsilon\right\}.
\end{eqnarray*}

Let $I$ be an interval of $\mathbb{R}$ and let $\{D(t)\}_{t\in I}$
be a family of subsets $D(t)$ of $\mathbb{R}^n$, we shall denote
$$D(I)=\bigcup_{t\in I}D(t)\times\{t\},$$
and $Q=D(\mathbb R)$.

Given a sufficiently smooth function $u$ of variable $(x,t)\in\mathbb{R}^{n+1}$,
we shall denote by
$\partial_i u=\frac{\partial u}{\partial x_i}$,
$\partial^2_{ij}u=\frac{\partial^2 u}{\partial x_i\partial x_j}$,
$i,j=1,\dots,n$ and $\partial_t u=\frac{\partial u}{\partial t}$.
For a multi-index $\beta=(\beta_1,\dots,\beta_n)$, $\beta_i\in\mathbb N$,
$i=1,\dots,n$ and $k\in\mathbb N$, we shall denote, as usual,
$\partial_x^{\beta}\partial_t^{k}u=\frac{\partial^{|\beta|+k}u}
{\partial x_1^{\beta_1}\dots\partial_{x_n}^{\beta_n}\partial_t^k}$,
where $|\beta|=\sum_{i=1}^n\beta_i$.
Also we shall write $\nabla=\nabla_x$, $\mathrm{div}=\mathrm{div}_x$.
For a matrix $A$ we shall denote by $A^*$ the transposed matrix of $A$.
We denote by
$\mathbb R^n_+=\{x\in\mathbb R^n\,:\,x_n>0\}$.

We shall use letters $C, C_0,C_1,\dots$ to denote constants.
The value of the constant can change from line to line,
but we shall specify the dependance everywhere they appear.
Sometimes we have dropped the dependance on $n$ which is fixed
($n\geq2$).

\bigskip

\noindent
\textbf{Functional Spaces}

\noindent
Let $D$ be a subset of $\mathbb R^{n+1}$, $f$ a function defined on $D$ with values
in $\mathbb R$ or $\mathbb R^n$ and $\alpha\in(0,1]$. We shall set
$$[f]_{\alpha;D}=\sup\left\{\frac{|f(x,t)-f(y,s)|}
{(|x-y|^2+|t-s|)^{\alpha/2}}
\,:\,(x,t),(y,s)\in D,\,(x,t)\neq(y,s)\right\}.$$
If $\alpha\in(0,2]$ we shall set
$$<f>_{\alpha;D}=\sup\left\{\frac{|f(x,t)-f(y,s)|}
{|t-s|^{\alpha/2}}
\,:\,(x,t),(y,s)\in D,\,t\neq s)\right\}.$$
Let $k$ be a positive integer number, $f$ a sufficiently smooth
function and $\alpha\in(0,1]$. We shall denote
$$[f]_{k+\alpha;D}=\sum_{|\beta|+2j=k}[\partial_x^{\beta}
\partial_t^j f]_{\alpha;D},\qquad
<f>_{k+\alpha;D}=\sum_{|\beta|+2j=k-1}[\partial_x^{\beta}
\partial_t^j f]_{1+\alpha;D}.$$

The following Sobolev spaces will be used
(we refer to \cite{Li-Ma} for further details).
We denote by $\Omega$ a bounded domain in $\mathbb R^n$.
The space $H=H^{3/2,3/4}_{,0}(\partial\Omega\times(0,T))$, its dual
$H'=H_1=H^{-3/2,-3/4}(\partial\Omega\times(0,T))$,
and $H_0=H^{1/2,1/4}(\partial\Omega\times(0,T))$.
We consider now the interpolation spaces between $H_0$ and $H_1$.
For any $\theta$, $0\leq\theta\leq1$, we define $H_{\theta}$ as $[H_0,H_1]_{\theta}$,
where the latter denotes
the interpolation at level $\theta$ between the two spaces $H_0$ and $H_1$.
The norm in $H_{\theta}$
will be denoted by
$\|\cdot\|_{\theta}$.
First, we notice that for any $\theta$, $0\leq\theta\leq1$,
there exists a constant $C_{\theta}$,
which depends on $\theta$ only, such that the following interpolation
inequality holds for any $\psi\in H_0$
\begin{equation}\label{interpineq}
\|\psi\|_{\theta}\leq C_{\theta}\|\psi\|^{1-\theta}_{0}\|\psi\|^{\theta}_{1}.
\end{equation}
We also make use of the following notation
$$W(\Omega\times(0,T))=\left\{v\,:\,
v\in L^2((0,T),H^1(\Omega)),\,\partial_t v\in L^2((0,T),H^{-1}(\Omega))
\right\}.$$

\bigskip

\noindent
\textbf{Boundary Regularity}

\noindent Let us give the following definitions

\begin{defin}
Let $\Omega$ be a domain in $\mathbb R^{n}$. Given
$\alpha$, $\alpha\in(0,1]$, we shall say
that $\partial\Omega$ is of class
$C^{1,\alpha}$ with constants $\rho_0,E>0$ if for any
$P\in\partial\Omega$, there exists a rigid transformation of $\mathbb R^n$
under which we have $P\equiv0$ and
$$\Omega\cap B_{\rho_0}(0)
=\{x\in B_{\rho_0}(0)\,:\,x_n>\varphi(x')\},$$
where $\varphi$ is a $C^{1,\alpha}$ function
on $B'_{\rho_0}(0)$ which satisfies the following conditions
$\varphi(0)=|\nabla_{x'}\varphi(0)|=0$ and
$\Vert\varphi\Vert_{C^{1,\alpha}(B'_{\rho_0}(0))}\leq E\rho_0$.
\end{defin}

\begin{oss}
We have chosen to normalize all norms in such a way that their terms are
dimensional homogeneous and coincide with the standard definition
when $\rho_0=1$.
For instance, for any $\varphi\in
C^{1,\alpha}(B'_{\rho_0}(0))$ we set
$$\Vert\varphi\Vert_{C^{1,\alpha}(B'_{\rho_0}(0))}
=\rho_0\Vert\varphi\Vert_{L^{\infty}(B'_{\rho_0}(0))}
+\rho_0 \Vert\nabla_{x'}\varphi\Vert_{L^{\infty}(B'_{\rho_0}(0))}
+\rho_0^{1+\alpha} [\nabla_{x'}\varphi]_{\alpha;B'_{\rho_0}(0)}$$
Similarly we shall set
$$\Vert u\Vert_{L^2(\Omega)}=\rho_0^{-(n+1)/2}\left(\int_{\Omega}u^2dX\right)^{1/2},$$
where $dX=dxdt$.
\end{oss}

\begin{defin}
Let $Q$ be a domain in $\mathbb R^{n+1}$. We shall say that
$Q$ (or equivalently $\partial Q$) is of class $\mathcal K$
with constants $\rho_0$, $E$ if for all $P_0\in\partial Q$
there exists a rigid transformation of space coordinates
under which we have $P_0=(0,0)$ such that
\begin{multline*}
Q\cap\left(B_{\rho_0}(0)\times(-\rho_0^2,\rho_0^2)\right)
=\{X\in B_{\rho_0}(0)\times(-\rho_0^2,\rho_0^2)\,:\,
x_n>\varphi(x',t)\},
\end{multline*}
where $\varphi$ is endowed with second derivatives with respect to $x_i$,
$i=1,\cdots,n$, with the $t$-derivative and with second derivatives
with respect to $x_i$ and $t$ and it satisfies the following conditions
$\varphi(0,0)=|\nabla_{x'}\varphi(0,0)|=0$ and
\begin{multline*}
\rho_0\Vert \nabla_{x'}\varphi\Vert_{L^\infty(B'_{\rho_0}\times(-\rho_0^2,\rho_0^2))}
+\rho_0^2\Vert D^2_{x'}\varphi\Vert_{L^\infty(B'_{\rho_0}\times(-\rho_0^2,\rho_0^2))}\\
+\rho_0^2\Vert \partial_t\varphi\Vert_{L^\infty(B'_{\rho_0}\times(-\rho_0^2,\rho_0^2))}
+\rho_0^3\Vert \nabla_{x'}\partial_t\varphi\Vert_{L^\infty(B'_{\rho_0}\times(-\rho_0^2,\rho_0^2))}
\leq E\rho_0.
\end{multline*}
\end{defin}

\begin{defin}[relative graphs]
\label{Defgrafrel}
Let $\beta\in(0,1]$.
We shall say that two bound\-ed
domains $\Omega _{1}$ and $\Omega _{2}$ in $\mathbb{R}^{n}$ of class $%
C^{1,\beta }$ with constants $R_{0},E$ are \textit{relative graphs} if for
any $P\in \partial \Omega _{1}$ there exists a rigid transformation of
coordinates under which we have $P\equiv 0$ and there exist $\varphi
_{P,1},\varphi _{P,2}\in C^{1,\beta }\left( B_{r_{0}}^{\prime }\left(
0\right) \right) $, where $\dfrac{r_{0}}{R_{0}}\leq 1$ depends on $E$ and $%
\beta $ only, satisfying the following conditions

i) $\varphi _{P,1}\left( 0\right) =0$ , $\left\vert \varphi _{P,2}\left(
0\right) \right\vert \leq \dfrac{r_{0}}{2}$,

ii) $\left\Vert \varphi _{P,i}\right\Vert _{C^{1,\beta }\left(
B_{r_{0}}^{\prime }\left( 0\right) \right) }\leq ER_{0},\quad i=1,2$,

iii) $\Omega _{i}\cap B_{r_{0}}\left( 0\right) =\left\{ x\in B_{r_{0}}\left(
0\right) :x_{n}>\varphi _{P,i}\left( x^{\prime }\right) \right\} $, $i=1,2$.

\noindent We shall denote
\begin{equation}
\gamma \left( \Omega _{1},\Omega _{2}\right) =\sup\limits_{P\in \partial
\Omega _{1}}\left\Vert \varphi _{P,1}-\varphi _{P,2}\right\Vert _{L^{\infty
}\left( B_{r_{0}}^{\prime }\left( 0\right) \right) }\text{ .}  \label{4.304}
\end{equation}
\end{defin}

\bigskip

\noindent
\textbf{The Dirichlet--to--Neumann map}

For any $g\in H$, let $u\in W(\Omega\times(0,T))$ be the weak solution to
the initial-boundary value problem
\begin{subequations}\label{dirpbm}
\begin{align}
\label{dirpbma}
&\partial_tu-\mathrm{div}((1+(k-1)\chi_Q)\nabla u)=0 &\text{in }\Omega \times(0,T),\\[1mm]
&u(x,0)=0 &\textrm{for }x\in\overline{\Omega},\\[1mm]
&u(x,t)=g(x,t) &\text{for }(x,t)\in\partial\Omega\times (0,T),
\end{align}
\end{subequations}
where $\chi_Q$ is the characteristic function of the set $Q$.

\noindent Then, for any $g\in H$, we set
\begin{equation}\label{lambdadef}
\Lambda_{Q}g=\frac{\partial u}{\partial\nu}|_{\partial\Omega\times(0,T)},
\quad u\text{ solution to }\eqref{dirpbm}.
\end{equation}
We have that there exists a unique solution $u\in W(\Omega\times(0,T))$
to problem \eqref{dirpbm} \cite{Ev}.
In addition, by standard regularity theorems \cite{Li}, \cite{la}
and by trace theorem \cite[Chapter~4, Theorem~2.1]{Li-Ma},
we conclude that $\Lambda_Q g$ belongs to the space $H_0$ and that the operator
$\Lambda_Q\,:\,H\to H_0$ is bounded.
We can also consider $\Lambda_{Q}$ as a linear and bounded operator between $H$ and $H'=H_1$, by setting
\begin{equation}\label{weakform}
\langle\Lambda_{Q}g,\phi\rangle_{H',H}=
\langle\frac{\partial u}{\partial\nu}|_{\partial\Omega\times(0,T)},\phi\rangle_{H',H}=
\int_{\partial\Omega\times(0,T)}\frac{\partial u}{\partial\nu}\phi,\quad\text{for any }g,\phi\in H,
\end{equation}
where $u$ solves \eqref{dirpbm} and
$\langle\cdot,\cdot\rangle_{H',H}$ is the duality pairing between $H'$ and $H$.

Let us remark that the operator $\Lambda_Q$ is usually referred to as the
\emph{Dirichlet--to--Neumann}
map associated to the equation \eqref{dirpbma}.

\subsection{Assumptions and Statement of the Main Result}\label{mt}
%
\noindent
\textbf{Assumptions on the domain}

\noindent
Let $\rho_0, M, E$ be given positive numbers.
We assume that $\Omega$ is a bounded domain in $\mathbb R^n$
satisfying
\begin{subequations}
\label{omega}
\begin{equation}
\label{omega1}
|\Omega|\leq M\rho_0^n,
\end{equation}
where $|\Omega|$ denotes the Lebesgue measure of $\Omega$.
We also assume that
\begin{equation}
\label{omega2}
\partial\Omega\textrm{ is of class } C^{1,1}
\textrm{ with constants } \rho_0, E.
\end{equation}
\end{subequations}

\noindent
\textbf{A priori information on the inclusion}

\noindent
Denoting by $Q=\bigcup_{t\in\mathbb R} D(t)\times\{t\}$ ($Q=D((-\infty,+\infty))$),
we assume the following conditions
\begin{subequations}
\label{q}
\begin{equation}
\label{q1}
\partial Q \textrm{ is of class } \mathcal K
\textrm{ with constans }\rho_0, E,
\end{equation}
\begin{equation}
\label{q2}
\textrm{dist}(D(t),\partial\Omega)\geq\rho_0,\qquad \overline{D(t)}\subset\Omega,
 \quad\forall\,t\in[0,T],
\end{equation}
\begin{equation}
\label{q3}
\Omega\setminus\overline{D(t)} \textrm{ is connected }\forall\,t\in[0,T].
\end{equation}
\end{subequations}

\begin{oss}
\label{Re4.1}
Let $t$ be any number in $[ 0,T] $. Observe that
\eqref{q1} automatically implies a lower bound on the diameter of every
connected component of $D(t)$ and $\Omega\setminus\overline D(t)$. In addition,
combining \eqref{omega1} with \eqref{omega2}, we have
an upper bound of the diameter of $\Omega$ and thus of
$D(t) $. Note also that \eqref{q1} and \eqref{q2}
implicitly comprise an a priori upper bound on the number of connected
components of $D(t)$.
\end{oss}

\begin{oss}
\label{n=1}
For $n=1$, it is possible to obtain Theorem \ref{maintheorem}
replacing assumption \eqref{q} by considering
$\Omega=(0,L)$ and
$D(t)=\{x\in\mathbb R\,:\, s_1(t)<x<s_2(t)\}$, where $s_i$, $i=1,2$, are $C^1$
functions such that for all $t\in(0,T)$
$$L-s_2(t)\geq\rho_0,\quad s_1(t)\geq\rho_0,\quad s_2(t)-s_1(t)\geq\rho_0$$
and
$$\|s_i\|_{L^\infty((0,T))}+\rho_0^2\|s'_i\|_{L^\infty((0,T))}\leq E\rho_0,
\qquad i=1,2,$$
here $3\rho_0<L$.
\end{oss}

\begin{teo}
\label{maintheorem}
Let $\Omega\subset\mathbb R^n$ satisfying \eqref{omega}.
Let $k>0$, $k\neq1$ be given.
Let $\{D_1(t)\}_{t\in\mathbb R}$, $\{D_2(t)\}_{t\in\mathbb R}$ be
two families of domains satisfying \eqref{q}.
Assume that for $\varepsilon>0$,
\begin{equation}
\label{ipoteo}
\|\Lambda_{Q_1}-\Lambda_{Q_2}\|_{\mathcal L(H,H')}\leq\varepsilon,
\end{equation}
where $Q_i=D_i((-\infty,+\infty))$, $i=1,2$.
Then
\begin{equation}
\label{tesiteo}
d_{\mathcal H}(\overline{D_1(t)},\overline{D_2(t)})\leq\omega_t(\varepsilon), \qquad t\in(0,T],
\end{equation}
where $\omega_t(s)$ is such that
\begin{equation}
\label{modulo}
\omega_t(s)\leq C\rho_0|\log s|^{-\eta},\qquad 0<s<1,
\end{equation}
with $C=C(t)$, depending on $t,M,E,k$ only,
and $0<\eta\leq1$ depending on $M,E,k$ only.
In addition we have that $C(t)$ tends to $+\infty$
as $t$ tends to $0$.
\end{teo}
Here $d_{\mathcal H}$ denotes the Hausdorff distance.

\section{Proof of Theorem \ref{maintheorem}}

For the sake of brevity we name $a_j=1+(k-1)\chi_{Q_j}$, $j=1,2$.
We fix $g\in H$.
We shall denote by $u_j$,
$j=1,2$ the solution of \eqref{dirpbm} when $Q=Q_j$.
For $\psi\in H^{1,1}(\Omega\times(0,T))$ such that
\begin{equation}
\label{psiT0}
\psi(\cdot,T)=0\quad\textrm{in }\Omega,
\end{equation}
using the weak formulation of \eqref{dirpbm} we have
\begin{multline*}
\int_{\partial\Omega\times(0,T)}a_j\frac{\partial u_j}{\partial\nu}\psi dS
+\int_{\Omega}u_j(x,0)\psi(x,0)dx\\
-\int_{\Omega\times(0,T)}\left(a_j\nabla u_j\cdot\nabla\psi
-u_j\partial_t\psi\right)dxdt=0\qquad\textrm{for }j=1,2.
\end{multline*}
Subtracting the two equations we obtain
\begin{multline}
\label{iniz2}
\int_{\Omega\times(0,T)}\left(a_1\nabla(u_1-u_2)\cdot\nabla\psi-(u_1-u_2)\partial_t\psi\right)dxdt\\
+\int_{\Omega\times(0,T)}(a_1-a_2)\nabla u_2\cdot\nabla\psi=
<(\Lambda_{Q_1}-\Lambda_{Q_2})g,\psi>_{H',H},
\end{multline}
(we notice here that in these identities it is possible to have
$u_i(\cdot,0)\neq0$ for $i=1,2$).
Taking $\psi$ such that it satisfies \eqref{psiT0} and
\begin{equation}
\label{psieq}
\partial_t\psi+\mathrm{div}(a_1\nabla\psi)=0
\qquad\textrm{ in }\Omega\times(0,T),
\end{equation}
by \eqref{iniz2} we have (recalling that on $\partial\Omega\times(0,T)$
$u_1=u_2=g$)
\begin{equation*}
\int_{\Omega\times(0,T)}(a_1-a_2)\nabla u_2\cdot\nabla \psi
=<(\Lambda_{Q_1}-\Lambda_{Q_2})g,\psi>_{H',H},\qquad\forall\,g\in H
\end{equation*}
or, equivalently,
\begin{equation}
\label{iniz4}
\int_0^T\int_\Omega(\chi_{Q_1}-\chi_{Q_2})
\nabla u_2\cdot\nabla\psi dxdt
=\frac{1}{k-1}<(\Lambda_{Q_1}-\Lambda_{Q_2})u_2,\psi>_{H',H}.
\end{equation}

Let us denote by $\Gamma_2(x,t;y,s)$ and $\Gamma_1^*(x,t;y,s)$
the fundamental solutions of the operator
$\partial_t-\mathrm{div}(a_2\nabla)$
and $\partial_t+\mathrm{div}(a_1\nabla)$
respectively
($\Gamma_1^*(x,t;y,s)=0$ if $t\geq s$ and $\Gamma_2(x,t;y,s)=0$ if $t\leq s$), that is
\begin{eqnarray*}
&&\int_{\mathbb R^{n+1}}\!\!\!\!\!\!\!\!
\left[-\Gamma_2(x,t;y,s)\partial_t\phi(x,t)
+a_2\nabla_x\Gamma_2(x,t;y,s)\cdot\nabla_x\phi(x,t)\right]dxdt
=\phi(y,s),\\[2mm]
&&\int_{\mathbb R^{n+1}}\!\!\!\!\!\!\!\!
\left[\Gamma_1^*(x,t;y,s)\partial_t\phi(x,t)
+a_1\nabla_x\Gamma_1^*(x,t;y,s)\cdot\nabla_x\phi(x,t)\right]dxdt
=\phi(y,s),
\end{eqnarray*}
for every $\phi\in C^1_0(\mathbb R^{n+1})$,
that is using the $\delta$ Dirac symbol, we have respectively
$$\partial_t \Gamma_2(x,t;y,s)-\mathrm{div}(a_2\nabla_x\Gamma_2(x,t;y,s))
=\delta(x-y,t-s)$$
and
$$-\partial_t \Gamma_1^*(x,t;y,s)-\mathrm{div}(a_1\nabla_x \Gamma_1^*(x,t;y,s))
=\delta(x-y,t-s).$$
Choosing in \eqref{iniz4} $u_2(x,t)=\Gamma_2(x,t;y,s)$
and $\psi(x,t)=\Gamma_1^*(x,t;\xi,\tau)$,
with $(y,s)$ and $(\xi,\tau)\notin \Omega\times(0,T)$,
$0\leq s<\tau\leq T$, we obtain
\begin{multline}
\label{iniz5}
\int_0^T\int_\Omega(\chi_{Q_1}-\chi_{Q_2})
\nabla_x \Gamma_2(x,t;y,s)\cdot\nabla_x \Gamma_1^*(x,t;\xi,\tau) dxdt\\
=\frac{1}{k-1}<(\Lambda_{Q_1}-\Lambda_{Q_2})
\Gamma_2(\cdot,\cdot;y,s),\Gamma_1^*(\cdot,\cdot;\xi,\tau)>_{H',H}.
\end{multline}
For $t\in[0,T]$ we shall define $\mathcal G(t)$ as the connected component
of $\Omega\setminus(\overline{D_1(t)}\cup\overline{D_2(t)})$ that contains
$\partial\Omega$, $\tilde{\mathcal G}(t)=(\mathbb R^n\setminus\Omega)\cup\mathcal G(t)$
and
$\tilde{\mathcal G}((0,T)):=\bigcup_{t\in(0,T)}\tilde{\mathcal G}(t)\times\{t\}$.
For $(y,s)$, $(\xi,\tau)\in\tilde{\mathcal G}((0,T))$ with $0\leq s<\tau\leq T$,
we set
\begin{eqnarray*}
&&S_1(y,s;\xi,\tau):=\int_{Q_1}\nabla_x\Gamma_2(x,t;y,s)
\cdot\nabla_x\Gamma_1^*(x,t;\xi,\tau)dxdt,\\[2mm]
&&S_2(y,s;\xi,\tau):=\int_{Q_2}\nabla_x\Gamma_2(x,t;y,s)
\cdot\nabla_x\Gamma_1^*(x,t;\xi,\tau)dxdt\\[2mm]
&&\mathcal U(y,s;\xi,\tau):= S_1(y,s;\xi,\tau)-S_2(y,s;\xi,\tau).
\end{eqnarray*}
\begin{oss}
\label{00}
Let us observe here that for $\tau<s$,
$S_1$ and $S_2$ can be defined identically
zero since
for $(y,s)$ and $(\xi,\tau)\in\tilde{\mathcal G}((0,T))$
we have
$S_1(y,s;\xi,s)=S_2(y,s;\xi,s)=0$.
\end{oss}
By \eqref{iniz5} we have
\begin{equation}
\label{iniz6}
\mathcal U(y,s;\xi,\tau)
=\frac{1}{k-1}<(\Lambda_{Q_1}-\Lambda_{Q_2})
\Gamma_2(\cdot,\cdot;y,s),\Gamma_1^*(\cdot,\cdot;\xi,\tau)>_{H',H},
\end{equation}
for all $y,\xi\notin\Omega$, $0\leq s<\tau\leq T$.
Denoting by $\Omega_D(t):=\Omega\setminus\overline{{\mathcal G}(t)}$,
$t\in[0,T]$, we introduce a variation of the Hausdorff distance
that, even though it is not a metric, we call it \emph{modified distance}
\begin{multline}
\label{dmu}
d_{\mu}(t)=d_\mu(D_1(t),D_2(t))\\
=\max\left\{\sup_{x\in\partial D_1(t)\cap\partial\Omega_D(t)}
\mathrm{dist}(x,D_2(t)),
\sup_{x\in\partial D_2(t)\cap\partial\Omega_D(t)}
\mathrm{dist}(x,D_1(t))\right\},
\end{multline}
$t\in[0,T]$.
We point out here that trivially
$d_\mu(D_1(t),D_2(t))\leq d_{\mathcal H}(\overline{D_1(t)},\overline{D_2(t)})$.
The following proposition shows the relation between the Hausdorff
distance and $d_\mu$, provided the priori assumptions \eqref{q} hold.
We refer to \cite[Proposition 3.3]{Al-DC} for the proof.
\begin{prop}
\label{dh}
Let $D_1(t)$ and $D_2(t)$ be two sets satisfying \eqref{q}
then for any $t\in[0,T]$
\begin{equation}
\label{dm}
d_{\mathcal H}(\partial D_1(t),\partial D_2(t))\leq C d_\mu(t),
\end{equation}
where $C$ depends on $E$ and $M$ only.
\end{prop}
We now give a proposition which connects the Hausdorff distance
between the boundary of the inclusions and between the inclusions.
\begin{prop}
\label{dm2}
Let $D_1$ and $D_2$ be two domains of class $C^2$ with
constants $E, \rho_0$ such that
$\mathbb R^n\setminus D_j$, $j=1,2$, is connected.
There exists a positive constant $C$ depending on $E$ only such that
\begin{equation}
\label{dm4}
d_{\mathcal H}(\overline D_1,\overline D_2)\leq C
d_{\mathcal H}(\partial D_1,\partial D_2).
\end{equation}
\end{prop}
{\bf Proof.} See Section \ref{distanza}

\begin{oss}
\label{rem8}
By Propositions \ref{dh} and \ref{dm2} we have
\begin{multline*}
d_{\mathcal H}(\overline D_1(t),\overline D_2(t))\leq
C d_{\mathcal H}(\partial D_1(t),\partial D_2(t))\leq\\
C'd_\mu(D_1(t),D_2(t))\leq
C'd_{\mathcal H}(\overline D_1(t),\overline D_2(t)).
\end{multline*}
Thus it turns out that the distances
 $d_{\mathcal H}(\overline D_1(t),\overline D_2(t))$,
$d_{\mathcal H}(\partial D_1(t),\partial D_2(t))$ and
$d_\mu(D_1(t),D_2(t))$ are equivalent.
\end{oss}

\begin{prop}
\label{grafrel}
Let $\Omega _{1}$ and $\Omega _{2}$ be bounded domains in $%
\mathbb{R}^{n}$ of class $C^{1,\beta }$ with constants $R_{0}$, $E$ and
satisfying $\left\vert \Omega _{i}\right\vert \leq MR_{0}^{n}$. There exist
numbers $\overline d$, $\overline\rho\in \left( 0,R_{0}\right] $ such that
$\dfrac{\overline d}{R_{0}}$ and $\dfrac{\overline\rho}{R_{0}}$ depend on $\beta $ and $E$
only, such that if we have
\begin{equation}
d_{\mathcal{H}}\left( \overline{\Omega _{1}},\overline{\Omega _{2}}\right)
\leq\overline d\text{,}  \label{4.360}
\end{equation}
then the following facts hold true

\noindent i) $\Omega _{1}$ and $\Omega _{2}$ are relative graphs and
\begin{equation}
\gamma \left( \Omega _{1},\Omega _{2}\right) \leq Cd_{\mathcal{H}}\left(
\overline{\Omega _{1}},\overline{\Omega _{2}}\right) \text{ ,}  \label{4.361}
\end{equation}
where $C$ depends on $\beta $ and $E$ only,

\noindent iii) any connected component of $\Omega _{1}\cap \Omega _{2}$ has boundary
of Lipschitz class with constants $\rho _{0}$, $L$, where $\rho _{0}$ is as
above and $L>0$ depends on $E$ only.
\end{prop}
{\bf Proof.} See \cite[Proposition 4.1.8]{Ve}.
\cvd

A key ingredient for the proof of our stability theorem are fundamental solutions.
We collect here some results we need.

We shall denote by $\Gamma_0(x-y,t-s)$ the standard fundamental solution of
$\partial_t-\Delta$ which is
$$\Gamma_0(x-y,t-s)=\frac{1}{[4\pi(t-s)]^{n/2}}
\mathrm e^{-\frac{|x-y|^2}{4(t-s)}},\qquad t>s.$$
We shall denote by $\Gamma(x,t;y,s)$ the fundamental solution
of the operator $\partial_t-\mathrm{div}((1+(k-1)\chi_Q)\nabla_x)$
(see \cite{ar}).
We recall that $\Gamma$ satisfies the following properties
\begin{equation}
\label{av3.1}
\Gamma(x,t;y,s)=\Gamma^*(y,s;x,t)\qquad\forall\,(x,t),(y,s)
\in Q,\,(x,t)\neq(y,s),
\end{equation}
where $\Gamma^*$ is the fundamental solution to $-\partial_t-\mathrm{div}((1+(k-1)\chi_Q)\nabla_x)$,
and
\begin{equation}
\label{av3.3}
0<\Gamma(x,t;y,s)\leq
\frac{C}{[4\pi(t-s)]^{n/2}}
\mathrm e^{-\frac{|x-y|^2}{C(t-s)}}\chi_{[s,+\infty)}(t),
\end{equation}
where $C\geq1$ depends on $k$ only.
Furthermore we have also the following estimate for the gradient of $\Gamma$.
\begin{prop}
\label{5-16Tpr}
Let $\Gamma(x,t;y,s)$ be the fundamental solution of the operator
$\partial_t-\mathrm{div}\left((1+(k-1)\chi_Q)\nabla_x\right)$.
There exists $C\geq1$, depending on $k$ and $E$ only such that
\begin{equation}
\label{5-16T}
\left|\nabla_x\Gamma(x,t;y,s) \right|\leq \frac{C}{(t-s)^{\frac{n+1}{2}}}
\,\mathrm e^{-\frac{|x-y|^2}{C(t-s)}},
\end{equation}
for almost every $x,y\in\mathbb R^n$ and $t,s\in\mathbb R$, $t>s$.
\end{prop}
{\bf Proof.} See Section \ref{fax-sec}.
\cvd

In the sequel we need the fundamental solution of the operator
$\mathcal L_+=\partial_t-\mathrm{div}((1+(k-1)\chi_+)\nabla)$
where $\chi_+=\chi_{\{(x,t)\in\mathbb R^{n+1}\,:\,x_n>0\}}$.
We shall denote by $\Gamma_+$ such a fundamental solution.
Also, we shall denote by $\Gamma_+^*$ the fundamental solution
of the adjoint operator of $\mathcal L_+$.
Observe that $\Gamma_+(x,t;y,s)=\Gamma_+(x,t-s;y,0)$
and $\Gamma_+^*(x,t;y,s)=\Gamma_+(x,s-t;y,0)$.
Here and in the sequel, for a given function $f(x',x_n)$,
we shall denote by $\mathcal F_{\zeta'}(f(\cdot,x_n))$ the Fourier transform
of $f$ with respect to the variable $x'$. Thus
$$\mathcal F_{\zeta'}(f(\cdot,x_n))=\int_{\mathbb R^{n-1}}f(x',x_n)
\mathrm e^{-ix'\cdot\zeta'}dx',$$
for every $\zeta'\in\mathbb R^{n-1}$.

In \cite{IKN} it has been proved some formulae for $\mathcal{F}_{\zeta
^{\prime }}\left( \Gamma _{+}\left( .,x_{n},t;y\right) \right) $. The
technique to prove such formulae is rather classical and lengthy. For this
reason we display only the formulae that we need corresponding to the case
in which $x_{n}>0$, $y_{n}<0$.

\underline{Case $k>1$}.

Denote by
\begin{eqnarray}
\label{(1)-2}
&&E(\zeta',x_n,t;\rho)=\exp\left[
-t(k-(k-1)\rho)|\zeta'|^2-\sqrt{\frac{k-1}{k}}x_n|\zeta'|\sqrt\rho\right],\\[2mm]
\label{(2)-2}
&&F(\zeta',y_n;\rho)=
\mathrm{Im}\left(A_1(\rho)\mathrm e^{iy_n\sqrt{k-1}\sqrt{1-\rho}|\zeta'|}\right),
\end{eqnarray}
where, for complex number $z=a+ib$, $\mathrm{Im}(z)$ denotes the imaginary part $b$ of $z$, and
\begin{equation}
\label{(3)-2}
A_1(\rho)=\frac{\sqrt{k-1}}{\pi}\frac{1}{i\sqrt{k-1}\sqrt{1-\rho}+\sqrt k\sqrt\rho}.
\end{equation}
Then
\begin{equation}
\label{(4)-2}
\mathcal F_{\zeta'}(\Gamma_+(\cdot,x_n,t;y,0))=
\int_0^1|\zeta'|\mathrm e^{-iy'\cdot\zeta'}E(\zeta',x_n,t;\rho)
F(\zeta',y_n;\rho)d\rho,
\end{equation}
for every $x_n>0$, $y_n<0$.

\underline{Case $0<k<1$}.

Denote by
\begin{eqnarray*}
&&G(\zeta',y_n,t;\rho)=
\exp\left[-t(1-(1-k)\rho)|\zeta'|^2+\sqrt{1-k}\,y_n|\zeta'|\sqrt\rho\right],\\[2mm]
&&H(\zeta',x_n;\rho)=\mathrm{Im}\left(A_2(\rho)\mathrm e^{-ix_n\sqrt{\frac{1-k}{k}}\sqrt{1-\rho}|\zeta'|}\right),
\end{eqnarray*}
where
$$A_2(\rho)=\frac{\sqrt{1-k}}{\pi}\frac{1}{\sqrt k\sqrt\rho-i\sqrt{1-k}\sqrt{1-\rho}}.$$
Then
$$\mathcal F_{\zeta'}(\Gamma_+(\cdot,x_n,t;y,0))=
\int_0^1|\zeta'|\mathrm e^{-iy'\cdot\zeta'}
G(\zeta',y_n,t;\rho)H(\zeta',x_n;\rho)d\rho,$$
for every $x_n>0$, $y_n<0$.
\begin{prop}
\label{pr1}
For every $\lambda_0\in(0,1]$ there exist
$\lambda_1,\lambda_2,\lambda_3\in(0,\lambda_0]$ such that
for every $h>0$ the following inequality holds true
\begin{multline}
\label{(1)-4}
I^{(h)}:=\left|\int_0^{\lambda_2 h^2}dt\int_{\mathbb R^n_+}
\nabla_x\Gamma_+^*(x,t;-\lambda_1h e_n,\lambda_2 h^2)\right.\\
\cdot\nabla_x\Gamma_0(x,t;-\lambda_3he_n,0)dx\Bigg|
\geq\frac{1}{Ch^n},
\end{multline}
where $C$, $C\geq1$, depends on $\lambda_1,\lambda_2,\lambda_3$ and $k$ only.
\end{prop}
{\bf Proof.} See Section \ref{subsec-pro-prefb}.\cvd

Through the paper we shall fix the value of $\lambda_1,\lambda_2,\lambda_3$ in such a way
that \eqref{(1)-4} is satisfied and we shall omit the dependence of various constants
by $\lambda_1,\lambda_2,\lambda_3$.
In the following we shall often make use of this technical lemma whose proof can be
found in \cite[Lemma 3 pg. 15]{fr}.
\begin{lem}
\label{fr}
Let $\alpha,\beta<\frac{n}{2}+1$ and $a>0$. Then
\begin{eqnarray*}
&&\int_s^t\int_{\mathbb R^n}(t-\tau)^{-\alpha}
\mathrm e^{-\frac{a|x-\xi|^2}{4(t-\tau)}}(\tau-s)^{-\beta}
\mathrm e^{-\frac{a|\xi-y|^2}{4(t-\tau)}}d\xi d\tau\\[2mm]
&=&\frac{C}{a^{n/2}}(t-s)^{\frac{n}{2}+1-\alpha-\beta}
\mathrm e^{-\frac{a|x-\xi|^2}{4(t-s)}},
\qquad\forall\,x,y\in\mathbb R^n,\,s<t,
\end{eqnarray*}
where $C$ depends on $\alpha$, $\beta$ and $n$ only.
\end{lem}


For $\overline t\in(0,T]$ fixed, we can assume, without loosing generality, that
there exists $O\in\partial D_1(\overline t)\cap\partial\Omega_D(\overline t)$
(for the sake of brevity we assume that $O$ is the origin of $\mathbb R^n$)
such that
\begin{equation}
\label{dmuO}
d_\mu(\overline t)=\mathrm{dist}(O,D_2(\overline t)).
\end{equation}
Denote by
$$\rho=\min\{d_\mu(\overline t),\rho_0\}.$$
Furthermore, denote by $\nu(O,\overline t)$ the exterior unit normal to $\partial D_1(\overline t)$ in $O$
pointing towards $\mathcal G(\overline t)$.
Now we introduce parameter $\delta\in(0,1]$
that we shall choose later on. We set
\begin{equation}
\label{t1y1yb}
t_1=\overline t-\lambda_2 h^2,\qquad \overline y=\lambda_1h\nu(0,\overline t),
\qquad y_1=\lambda_3 h\nu(0,\overline t),
\end{equation}
where
\begin{equation}
\label{3.24bis}
0<h\leq\delta\min\{\rho,\sqrt{\overline t}\}.
\end{equation}
Notice that \eqref{3.24bis} implies that $t_1\in(0,\overline t)$.
By using \eqref{q1} it is simple to check that there exists $C_1$,
$C_1\geq1$, depending on $E$ only such that if
\begin{equation}
\label{1-14T}
0<\delta\leq\frac{\lambda_3}{C_1}
\end{equation}
then, for every $t\in[t_1,\overline t]$, we have
\begin{eqnarray}
\label{2a-14T}
&&\mathrm{dist}(\overline y,D_1(t))\geq\frac{1}{2}
\min\left\{\lambda_1,\lambda_2,\lambda_3\right\}h,\\[2mm]
\label{2b-14T}
&&\mathrm{dist}(y_1,D_1(t))\geq\frac{1}{2}
\min\left\{\lambda_1,\lambda_2,\lambda_3\right\}h.
\end{eqnarray}
On the other side, using the inequality
\cite[Proposition 4.1.6]{Ve}
\begin{equation}
\left|\mathrm{dist}(O,D_2(t))-\mathrm{dist}(O,D_2(\overline t))\right|
\leq\frac{C_0}{\rho_0}|t-\overline t|,
\end{equation}
where $C_0$ depends on $E$ and $M$ only, for $t\in[t_1,\overline t]$
and by using the triangle inequality
we have that there exists $C_2$, $C_2\geq1$, depending on $E$ and $M$ only
such that if
\begin{equation}
\label{(1)-12}
0<\delta\leq\frac{1}{C_2}
\end{equation}
then for $t\in[t_1,\overline t]$
\begin{equation}
\label{3-16T}
\mathrm{dist}(z,D_2(t))\geq\frac{1}{2}\rho, \qquad\textrm{with }z=\overline y,y_1.
\end{equation}
\begin{prop}
\label{prefb}
Let $\{D_1(t)\}_{t\in\mathbb R}$, $\{D_2(t)\}_{t\in\mathbb R}$ be two families of
domains satisfying \eqref{q}
and let $\lambda_1,\lambda_2,\lambda_3\in(0,1)$ be such that
the inequality \eqref{(1)-4} is satisfied.
Then there exist $C$, $C\geq1$, and $\tilde C$, $\tilde C\geq1$,
$C$ depending on $k$ only and $\tilde C$ depending on
$k,E,M,\lambda_1,\lambda_2$ and $\lambda_3$ only such that
\begin{equation}
\label{efb}
|\mathcal U(y_1,t_1;\overline y,\overline t)|\geq
\frac{1}{Ch^n},
\end{equation}
for $0<h\leq\frac{1}{\tilde C}\min\{\rho,\sqrt{\overline t}\}$,
where $y_1,t_1,\overline y,\overline t$, and $\rho$ are defined in \eqref{t1y1yb}.
\end{prop}
{\bf Proof.} See Section \ref{subsec-pro-prefb}
\cvd

\begin{teo}[Two-spheres and one-cylinder inequality]
\label{2sfe1cil}
Let $\lambda$, $\Lambda$ and $R$ positive numbers with $\lambda\in(0,1]$.
Let $P$ be the parabolic operator
\begin{equation}
\label{3.e1}
P=\partial_t-\partial_i\left(a^{ij}\partial_j\right),
\end{equation}
where $\{a^{ij}(x,t)\}_{i,j=1}^n$ is a symmetric $n\times n$ matrix.
For $\xi\in\mathbb R^n$ and $(x,t),(y,s)\in\mathbb R^{n+1}$ assume that
\begin{subequations}
\label{3.e23}
\begin{equation}
\label{3.e2}
\lambda|\xi|^2\leq\sum_{i,j=1}^na^{ij}(x,t)\xi_i\xi_j\leq\lambda^{-1}|\xi|^2
\end{equation}
and
\begin{equation}
\label{3.e3}
\left(\sum_{i,j=1}^n\left(a^{ij}(x,t)-a^{ij}(y,s)\right)^2\right)^{1/2}
\leq\frac{\Lambda}{R}\left(|x-y|^2+|t-s|\right)^{1/2}.
\end{equation}
\end{subequations}
Let $u$ be a function in $H^{2,1}\left(B_R\times(0,R^2)\right)$
satisfying the inequality
\begin{equation}
\label{3.e4}
|Pu|\leq\Lambda\left(\frac{|\nabla u|}{R}+\frac{|u|}{R^2}\right)
\qquad\textrm{in }B_R\times(0,R^2].
\end{equation}
Then there exist constants $\eta_1\in(0,1)$ and $C\in[1,+\infty)$,
depending on $\lambda$, $\Lambda$ and $n$ only such that
for every $r_1,r_2$, $0<r_1\leq r_2\leq\eta_1R$ we have
\begin{equation}
\label{3.e5}
\|u(\cdot,R^2)\|_{L^2(B_{r_2})}\leq
\frac{CR}{r_2}\|u\|^{1-\theta_1}_{L^2(B_R\times(0,R^2))}
\|u(\cdot,R^2)\|_{L^2(B_{r_1})}^{\theta_1},
\end{equation}
where $\theta_1=\frac{1}{C\log\frac{R}{r_1}}$.
\end{teo}
{\bf Proof.} See \cite{Ve}
\cvd

\bigskip

We can now start to prove our stability theorem. Before entering into details,
we wish to warn the reader that, sometimes we use the previous auxiliary results
(such as Lemma \ref{fr}  or Proposition \ref{5-16Tpr})
omitting some computations that are similar to the one contained in the proofs
of Section 4.

\bigskip

\noindent{\bf Proof of Theorem \ref{maintheorem}.}
We divide the proof of the theorem in two steps.
In the first step we provide a rough stability estimate (see \eqref{rough} below).
In the second step we prove the logarithmic stability estimate \eqref{modulo}.

\bigskip

\noindent{\bf Step 1.}

\noindent  We shall denote by
\begin{equation}
\label{xi}
\Xi_{\rho_0}=\{x\in\mathbb R^n\,:\,\rho_0/2<\mathrm{dist}(x,\Omega)<\rho_0\}
\end{equation}
and
\begin{equation}
\label{omega_ro}
\Omega_{\rho_0}=\{x\in\mathbb R^n\,:\,\mathrm{dist}(x,\Omega)<\rho_0\}.
\end{equation}
Since $\mathcal U(y,s;\xi,\tau)$ is equal to $0$ when $s\geq\tau$
(see Remark \ref{00}),
for $(y,s)\in\Xi_{\rho_0}\times(0,T)$ we define
\begin{equation}
\label{(1)-10}
v(\xi,\tau):=\mathcal U(y,s;\xi,\tau),\qquad
(\xi,\tau)\in\tilde{ \mathcal G}((0,T)):=\cup_{t\in(0,T)}\tilde{\mathcal G}(t)\times\{t\}.
\end{equation}
Let $h$ be the one defined in \eqref{3.24bis},
with $\delta\leq\frac{1}{\tilde C}$
and let $\lambda_1,\lambda_2,\lambda_3\in(0,1)$
be such that inequality \eqref{(1)-4} is satisfied.
Let $\overline x\in\Xi_{\rho_0}$ be such that
$\mathrm{dist}(\overline x,\mathbb R^n\setminus\Xi_{\rho_0})\geq\rho_0/8$.
Let us denote by $\gamma$ a simple connected arc in
$(\Omega_{\rho_0}\setminus\Omega_D(\overline t))_{\lambda_1 h/2}
=\{x\in\Omega_{\rho_0}\setminus\Omega_D(\overline t)\,:\,
\mathrm{dist}(x,\partial(\Omega_{\rho_0}\setminus\Omega_D(\overline t)))>\frac{\lambda_1h}{2}\}$,
connecting $\overline x$ to $\overline y$,
where $\overline y$ is defined in \eqref{t1y1yb}.
By \eqref{ipoteo} we have
\begin{equation}
\label{1)-2}
\| v(\cdot,\overline t)\|_{L^\infty(B_{\rho_0/2}(\overline x))}\leq\frac{C\varepsilon}{\rho_0^n},
\end{equation}
and by Lemma \ref{fr} and Proposition \ref{5-16Tpr} we have
\begin{equation}
\label{2)-2}
\|v\|_{L^{\infty}(\mathcal G((0,T)))}\leq\frac{C}{\rho_0^n},
\end{equation}
where $C$ depends on $k$ only.
It is easy to check that by \eqref{omega} and \eqref{q}
there exists $C$, $C\geq\tilde C$, depending on $k,E,M$
such that for all
$x\in(\Omega_{\rho_0}\setminus\Omega_D(\overline t))_{\lambda_1h/2}$
and $0<h\leq\frac{1}{C}\min\{\rho,\sqrt{\overline t}\}$,
\begin{equation}
\label{1)-3}
Q_{\lambda_1h/4}:=
B_{\lambda_1h/4}(x)\times\left(\overline t-\left(\frac{\lambda_1h}{4}\right)^2,\overline t\right]
\subset\tilde{\mathcal G}((-\infty,\overline t]).
\end{equation}
Since $v$ solves the heat equation, we can apply Theorem \ref{2sfe1cil}
along a chain of balls centered in points of $\gamma$.
More precisely, let us define $\overline\rho=\eta_1\lambda_1h/12$, where $\eta_1\in (0,1)$
is defined in Theorem \ref{2sfe1cil}, and $x_i$, $i=1,\dots,m_h$
as follows: $x_1=\overline x$, $x_{i+1}=\gamma(t_i)$,
where $t_i=\max\{t\,:\,|\gamma(t)-x_i|=2\overline\rho\}$,
if $|x_i-\overline y|>2\overline\rho$,
otherwise let $i=m_h$ and stop the process.
We have
$m_h\leq CM\left(\frac{\rho_0}{h}\right)^n$, where
$C>0$ is an absolute constant.
By construction the balls $B_{\overline\rho}(x_i)$ are pairwise disjoint
and $|x_{i+1}-x_i|=2\overline\rho$ for $i=1,\dots,m-1$ and
$|x_m-\overline x|\leq2\overline\rho$.
By an iterated application of the two--sphere and one--cylinder inequality (Theorem \ref{2sfe1cil})
to $v$ with $R=\lambda_1h/4$, $r_2=3\overline\rho$, $r_1=\overline\rho$
over the chain of balls $B_{\overline\rho}(x_i)$,
since we have $B_{r_1}(x_{i+1})\subset B_{r_2}(x_i)$, $i=1,\dots, m_h$,
by \eqref{1)-2} and \eqref{2)-2}
we have
\begin{equation}
\label{1s}
\left(\frac{1}{r_2^n}\int_{B_{r_2}(\overline y)}
v^2(\xi,\overline t)d\xi\right)^{1/2}
\leq  \frac{C}{\rho_0^n}\varepsilon^{s_2^{m_h}}(C+\varepsilon)^{1-s_2^{m_h}},
\end{equation}
where $s_2\in(0,1)$ is an absolute constant and $C$ depends on a priori data only.
From now on, in order to simplify the writing and since the case $\varepsilon\geq1$
is trivial, we shall assume that $\varepsilon\in(0,1)$.
By standard regularity estimates \cite{Li} and taking into account
\eqref{1)-2} and \eqref{2)-2} we have
\begin{equation}
\label{1)-4}
\|\nabla v(\cdot,\overline t)\|_{L^\infty(B_{r_2/2}(\overline y))}
\leq\frac{C}{h\rho_0^n},
\end{equation}
where $C$ depends on $k$ only.
Recalling now the interpolation inequality
(see \cite[(5.29)]{Al-Be-Ro-Ve})
\begin{equation}
\label{inter}
\|g\|_{L^\infty(B_r)}\leq
C\left[\|g\|_{L^\infty(B_r)}+
r\|\nabla g\|_{L^\infty(B_r)}\right]^{\frac{n}{n+2}}
\left(r^{-n}\int_{B_r}g^2\right)^{\frac{1}{n+2}},
\end{equation}
for every $r>0$,
where $C$ is an absolute constant,
by \eqref{1)-4} and \eqref{1s} we have
\begin{equation}
\label{2)-4}
\|v(\cdot,\overline t)\|_{L^\infty(B_{\frac{r_2}{2}}(\overline y))}
\leq\frac{C}{\rho_0^n}\varepsilon^{s_2^{m_h}}
:=\omega_h^{(1)}(\varepsilon),
\end{equation}
where $C$ depends on the a priori data only.
Now defining $w(y,s)=\mathcal U(y,s;\overline y,\overline t)$
and taking into account \eqref{2)-4} we have
$$\left\{\begin{array}{ll}
\partial_s w(y,s)+\Delta_y w(y,s)=0,&\textrm{ in }\tilde{\mathcal G}((0,T)),\\[2mm]
w(y,s)_{|\{s\geq\overline t\}}=0,\\[2mm]
|w(y,t)|\leq\omega_h^{(1)}(\varepsilon),&(y,s)\in\Xi_{\rho_0}\times(0,T).
\end{array}\right.$$
Now we want to estimate from above $|w(y_1,t_1)|$.
In order to obtain such an estimate we argue as before, but here,
instead of \eqref{2)-2}, we use the inequality
\begin{equation}
\label{1)-5}
\|w\|_{L^\infty(\tilde{\mathcal{ G}}^{(h)}([t_1,\overline t]))}\leq\frac{C}{h^n},
\end{equation}
where $\tilde{\mathcal G}^{(h)}([t_1,\overline t])=\{(x,t)\in\tilde{\mathcal G}((0,T))\,:\,
t_1\leq t\leq\overline t,\,\mathrm{dist}(x,\Omega_D(t))\geq\delta_2 h\}$,
$\delta_2=\frac{1}{8}\min\{\lambda_1,\lambda_2,\lambda_3\}$
and $C$ depends on $k$ only.
Inequality \eqref{1)-5} is a consequence of Proposition \ref{5-16Tpr}
and Lemma \ref{fr}. Notice that by virtue of \eqref{2a-14T} and \eqref{2b-14T}
we have $(y_1,t_1),(\overline y,\overline t)\in\tilde{\mathcal G}^{(h)}([t_1,\overline t])$.
Therefore we obtain
\begin{equation}
\label{u1}
|\mathcal U(y_1,t_1;\overline y,\overline t)|
=|w(y_1,t_1)|
\leq\frac{C}{h^n}\varepsilon^{\tilde s_2^{m_h}},
\end{equation}
where $\tilde s_2\in(0,1)$ is an absolute constant and $C$ depends
on the a priori data only.
Now we introduce some notation.
We set $\delta(\overline t)=\min\left\{\frac{\sqrt{\overline t}}{\rho_0},1\right\}$
and $h/\rho_0=q$. By Proposition \ref{prefb} and \eqref{u1} we have that there exists
$C_5$, $C_5\geq1$, depending on $k,E$ and $M$ only such that
\begin{equation}
\label{u3}
1\leq C_5\varepsilon^{s_3^{q^{-n}}},
\qquad\textrm{for every }
q\leq\frac{1}{C_5}\min\left\{\frac{d_\mu(\overline t)}{\rho_0},\delta(\overline t)\right\},
\end{equation}
where $s_3$, $s_3\in(0,1)$, depends on $M$ only.
We distinguish two cases
$$ \mathrm i) \,d_\mu(\overline t)\geq\min\{\sqrt{\overline t},\rho_0\}\qquad
\mathrm{ii})\,d_\mu(\overline t)<\min\{\sqrt{\overline t},\rho_0\}.$$
If case i) occurs we choose $q$ such that $s_3^{q^{-n}}=|\log\varepsilon|^{-1/2}$,
that is
$$q=q_\varepsilon:=\left(\frac{|\log s_3|}
{|\log |\log \varepsilon|^{-1/2}|}
\right)^{\frac{1}{n}}.$$
Denote by $\varepsilon _{\ast}(\overline t)$ the least upper
bound of the set
$\{\varepsilon \in (0,1)\,:\,q_{\varepsilon}
\leq \delta(\overline t)\}$.
By \eqref{u3} we have
$$1\leq C_{5}\exp \left\{-\left\vert \log \varepsilon \right\vert
^{1/2}\right\},$$
which, for $0<\varepsilon \leq \varepsilon _{\ast\ast}(t)
:=\min \left\{ \varepsilon _{\ast }\left( \overline{t}\right),
\mathrm e^{-\left( \log C_{5}\right) ^{2}}\right\} $,
yields to a contradiction.
Thus, if $0<\varepsilon \leq \varepsilon _{\ast\ast}(\overline t)$,
case i) cannot occur.

\noindent Let us consider now case ii), that is
$d_{\mu}(\overline t)<\min \left\{ \sqrt{\overline{t}},\rho _{0}\right\} $.
By (3.46) we have
$$1\leq C_{5}\exp \left\{ -s_{3}^{q^{-n}}\left\vert \log \varepsilon
\right\vert \right\},$$
for every $q\leq \dfrac{C_{5}^{-1}d_{\mu }\left( \overline{t}\right) }{\rho
_{0}}$. Now, if
\begin{equation}
d_{\mu }\left( \overline{t}\right) \leq 2C_{5}\rho _{0}\left\vert \log
\left( \left\vert \log \varepsilon \right\vert ^{\frac{\left\vert \log
s_{3}\right\vert ^{-1}}{2}}\right) \right\vert ^{-\frac{1}{n}}
\label{1}
\end{equation}
then we are done. On the other side, if
\begin{equation}
d_{\mu }\left( \overline{t}\right) >2C_{5}\rho _{0}\left\vert \log \left(
\left\vert \log \varepsilon \right\vert ^{\frac{\left\vert \log
s_{3}\right\vert ^{-1}}{2}}\right) \right\vert ^{-\frac{1}{n}},
\label{2}
\end{equation}
let us denote by
$$\widetilde{q}_{\varepsilon }=\left\vert \log \left( \left\vert \log
\varepsilon \right\vert ^{\frac{\left\vert \log s_{3}\right\vert ^{-1}}{2}
}\right) \right\vert ^{-\frac{1}{n}}$$
and 
by $\varepsilon _{0}\left( \overline{t}\right) $ the least
upper bound of the set $\left\{ \varepsilon \in \left( 0,\varepsilon _{\ast
}\left( \overline{t}\right) \right) :\widetilde{q}_{\varepsilon }\leq \delta
\left( \overline{t}\right) \right\} $. Now, for $0<\varepsilon \leq
\varepsilon _{0}\left( \overline{t}\right) $, we choose
$q=\widetilde{q}_{\varepsilon }$ and by \eqref{u3} we have
$$1\leq C_{5}\exp \left\{ -\left\vert \log \varepsilon \right\vert
^{1/2}\right\}.$$
Since the last inequality yields to a contradiction whenever
$0<\varepsilon\leq \varepsilon _{0}\left( \overline{t}\right) $, we have that if
$0<\varepsilon \leq \varepsilon _{0}\left( \overline{t}\right) $ then (\ref{2})
cannot occur, so inequality (\ref{1}) holds true . Finally, by using
Proposition \ref{dh} and \ref{dm2} and Remark \ref{rem8}, we have
\begin{equation}
\label{rough}
d_{\mathcal H}(\overline{D_1(\overline t)},\overline{D_2(\overline t)})\leq
2C_5\rho_0
\left|\log\left(|\log\varepsilon|^{\frac{|\log s_3|^{-1}}{2}}\right)
\right|^{-1/n}:=\sigma(\varepsilon),
\end{equation}
for $0<\varepsilon\leq\varepsilon_0(\overline t)$.

\bigskip

\noindent{\bf Step 2.}

In order to prove estimate \eqref{tesiteo} we apply Proposition \ref{grafrel}
to $\Omega_i:=\Omega\setminus\overline{D_i(\overline t)}$, $i=1,2$,
and $R_0=\rho_0$.
Indeed by \eqref{rough} we have that, for $\varepsilon$ small enough,
$\Omega_1$ and $\Omega_2$ are relative graphs. More precisely if
$0<\varepsilon\leq\min\{\varepsilon_0(\overline t),\overline d\}$
(where $\overline d$ is defined in Proposition \ref{grafrel})
then there exists
$r_0>0$
($r_0/\rho_0\leq1$ depending on $E$ only) such that
\begin{equation}
\label{ri1}
B_{r_0}(0)\cap D_i(\overline t)
=\{x\in B_{r_0}(0)\,:\, x_n>\varphi_i(x')\},
\qquad i=1,2,
\end{equation}
and $\|\varphi_1-\varphi_2\|_{L^\infty(B'_{r_0})}
\leq C\rho_0\sigma(\varepsilon)$,
where $C$ depends on $E$ only. By \eqref{q1}
and an interpolation inequality \cite[(5.30)]{Al-Be-Ro-Ve} we have that
$\|\varphi_1-\varphi_2\|_{C^{1}(B'_{r_0})}
\leq C\rho_0(\sigma(\varepsilon))^{\overline\beta}$,
$\overline\beta\in(0,1)$.
Thus, with eventually a rigid transform of coordinates,
provided we pick a smaller $r_0$, there exists
$\varepsilon_0>0$, depending on $E$ only,
such that for $\varepsilon\leq\varepsilon_0$
we can assume that $|\nabla\varphi_1(0)|=0$.
In the sequel we continue to denote by $\varepsilon_0(\overline t)$
the number $\min\{\varepsilon_0(\overline t),\varepsilon_0\}$.
Let us define, for a unit vector $\zeta$ and $0<\alpha<\pi/2$
$$\mathcal C(z,\zeta,\alpha,r_0)
=\left\{x\in B_{r_0}(z)\,:\, \frac{(x-z)\cdot\zeta}{|x-z|}>\cos\alpha \right\}.$$
By \eqref{ri1} we have that
$\mathcal C(0,\nu,\alpha,r_0)\subset\mathcal G(\overline t)$,
where $\alpha$, $\alpha\in(0,\pi/2)$, depends on $E$ only.
Let us denote $\overline\mu=\frac{\rho_0}{1+\sin\alpha}$,
$\delta^*=\frac{\sin\alpha}{\sqrt{2E+1}}$ and $\overline\rho_0=\overline\mu\cos\frac{\alpha}{2}$.
We have that $\mathcal S((0,\overline t),\nu,\frac{\alpha}{2},\delta^*,\overline\rho)
\subset\tilde{\mathcal G}((-\infty,\overline t])$, where we set $\nu=\nu(0,\overline t)$,
for the sake of brevity, and
\begin{multline*}
\mathcal S\left((0,\overline t),\nu,\frac{\alpha}{2},\delta^*,\overline\rho\right)\\
=\left\{(z,t)\in\mathbb R^{n+1}\,:\,
z\in\mathcal C(0,\nu,\frac{\alpha}{2},\overline\rho),
\,\overline t-\delta^*(x\cdot\nu)^2< t\leq\overline t \right\}
\end{multline*}
and $\nu=\nu(0,\overline t)$.
We want to estimate $v(\overline y,\overline t)=\mathcal U(y,s;\overline y,\overline t)$
when $(y,s)\in\Xi_{\rho_0}\times(0,T),$ $s<\overline t$, where $v$ solves
\begin{equation}
\label{v}
\left\{\begin{array}{ll}
\partial_\tau v-\Delta_\xi v=0,&\textrm{for }(\xi,\tau)\in\tilde{\mathcal G}((0,T)),\\[2mm]
v(\xi,\tau)_{|\tau\leq s}=0,\\[2mm]
\|v\|_{L^\infty(\Xi_{\rho_0}\times(0,T))}\leq\tilde\varepsilon,
\end{array}\right.
\end{equation}
where $\tilde\varepsilon=C\varepsilon/\rho_0^n$,
where $C$ depends on $k$ only.
Let us fix $(y,s)$ such that
\begin{equation}
\label{star}
(y,s)\in\tilde{\Xi}_{\rho_0}\times(0,T)=
\{x\in\mathbb R^n\,:\,\frac{5}{8}\rho_0
\leq d(x,\Omega)\leq\frac{7}{8}\rho_0\}\times(0,T).
\end{equation}
By Proposition \eqref{5-16Tpr} and Lemma \eqref{fr},
taking into account the last relation of \eqref{v} and by \eqref{star}
we have
\begin{equation}
\label{1-222}
\|v\|_{L^\infty(\tilde{\mathcal G}(0,T))}\leq\frac{C}{\rho_0^n}+\tilde\varepsilon:=H.
\end{equation}
In the sequel we continue to denote by $v$ the trivial extension of $v$.
Notice that, by \eqref{1-222}, we have
$$\| v\|_{L^\infty(\mathcal S((0,\overline t),\nu,\alpha/2,\delta^*,\overline\rho)}\leq H.$$
Denote by
\begin{equation*}
\alpha _{1}=\arcsin \left( \min \left\{ \sin \frac{\alpha }{2},\delta^*
\left( 1-\sin \frac{\alpha }{2}\right) \right\} \right),
\end{equation*}
\begin{equation*}
\mu _{1}=\frac{\overline{\rho }}{1+\sin \alpha _{1}},w_{1}=z+\mu
_{1}\zeta,\rho _{1}=\frac{1}{4}\mu _{1}\eta _{1}\sin \alpha _{1},
\end{equation*}
where $\eta _{1}\in \left( 0,1\right) $ is defined in Theorem \ref{2sfe1cil}.
We have
\begin{equation}
\label{4.205}
\mathrm{dist}\left( w_{1},\partial\mathcal G\left( \overline t\right) \right) \geq \min
\left\{\rho _{0}-\left\vert w_{1}-z\right\vert ,\left\vert w_{1}-z\right\vert \sin
\alpha \right\} =\rho _{0}\widetilde{\eta },
\end{equation}
where
\begin{equation*}
\widetilde{\eta }=\min \left\{ 1-\frac{\cos \frac{\alpha }{2}}{1+\sin \alpha
}\frac{1}{1+\sin \alpha _{1}},\frac{\sin \alpha }{1+\sin \alpha _{1}}
\frac{\cos \frac{\alpha }{2}}{1+\sin \alpha }\right\}.
\end{equation*}
Now $\left( \mathcal G\left( \overline t\right) \right) _{\frac{\rho _{0}\widetilde{\eta }}{2}}$
is connected and, by (\ref{4.205}), $w_{1}\in \left(\mathcal G\left(
\overline t\right) \right) _{\frac{\rho _{0}\widetilde{\eta }}{2}}$.
Therefore by an iterated application of the two-sphere and one-cylinder inequality
(see also \cite[Proposition 4.1.1]{Ve}) we get
\begin{equation}
\label{4.210}
\left( \rho _{1}^{-n}\int\nolimits_{B_{\rho _{1}(w_1)}}
v^{2}\left( \xi,\overline t\right) d\xi\right)^{1/2}\leq
C\tilde\varepsilon^{s_{4}}H^{1-s_{4}},
\end{equation}
where $s_{4}$, $s_{4}\in \left( 0,1\right)$, is an absolute constant
and $C$ depends on the a priori data only.
Denote
\begin{eqnarray*}
&&\mu_k=a^{k-1}\mu_1,\quad w_k=\mu_k\nu \quad \rho_k=a^{k-1}\rho_1\\
&&d_k=\mu_k-\rho_k=a^{k-1}\mu_1(1-\frac{1}{4}\eta_1\sin\alpha_1),
\end{eqnarray*}
where $a=\frac{1-\frac{1}{4}\eta_1\sin\alpha_1}
{1+\frac{1}{4}\eta_1\sin\alpha_1}$.
For every $k\geq1$, the following inclusions hold true
\begin{eqnarray}
\label{1b}
&&B_{\rho_{k+1}}(w_{k+1})\subset B_{3\rho_k}(w_k)\subset B_{4\eta_1^{-1}\rho_k}(w_k)
\subset\mathcal C(0,\nu,\alpha_1,r_0)\\
\label{2b}
&&B_{4\eta_1^{-1}}(w_k)\times(\overline t-(4\eta_1^{-1}\rho_k)^2,\overline t]
\subset\mathcal S((0,\overline t),\nu,\alpha,\delta^*,r_0).
\end{eqnarray}
Let us consider $h$ defined in \eqref{3.24bis}.
We further assume $\lambda_1h\in(0,d_1]$.
Let $\overline{k}$ be the smallest positive integer such that $d_{k}\leq\lambda_1 h$. We
have
\begin{equation}
\label{3.870}
\frac{\left\vert \log \left(\lambda_1 h/d_{1}\right) \right\vert }{\left\vert \log
a\right\vert }\leq \overline{k}-1\leq \frac{\left\vert \log \left(\lambda_1
h/d_{1}\right) \right\vert }{\left\vert \log a\right\vert }+1\text.
\end{equation}
Denote
\begin{equation*}
\sigma _{j}=\left( \rho_0^{-n}
\int\nolimits_{B_{\rho _{j}}\left( w_{j}\right)}
v^{2}\left(\xi,\overline t\right) d\xi\right) ^{1/2}, \quad j=1,...,\overline{k}.
\end{equation*}
By the Theorem \ref{2sfe1cil},
\eqref{1b}, \eqref{2b} and since
\begin{equation*}
\sigma _{j+1}\leq\left( \rho_0^{-n}
\int\nolimits_{B_{3\rho _{j}}\left( w_{j}\right)}
v^{2}\left(\xi,\overline t\right) d\xi\right) ^{1/2},\quad j=1,...,\overline{k}-1,
\end{equation*}
we obtain
\begin{equation}
\label{1-16}
\sigma^2_{j+1}\leq C_1H^{2(1-\theta_*)}\sigma_j^{2\theta_*},\quad j=1,\dots,\overline k-1,
\end{equation}
where $\theta_*=\frac{1}{C_0\log(4/\eta_1)}$.
By iterating \eqref{1-16} we get
\begin{equation}
\label{1)-10}
\sigma_{\overline k}^2\leq C^{\frac{1}{1-\theta_*}}H^{2(1-\theta_*^{\overline k})}
\sigma_1^{2\theta_*^{\overline k}}.
\end{equation}
By Lemma \ref{fr} and Proposition \ref{5-16Tpr} we have
\begin{equation}
\label{1)-11}
\|v\|_{L^\infty(Q_{\rho_{\overline k}})}\leq\frac{C}{\rho_0^n},
\end{equation}
where $Q_{\rho_{\overline k}}=B_{\rho_{\overline k}}(w_{\overline k})
\times(\overline t-\rho_{\overline k}^2,\overline t]$.
By standard regularity estimate and \eqref{1)-11} we get
\begin{equation}
\label{2)-11}
\|\nabla v(\cdot,\overline t)\|_{L^\infty(B_{\rho_{\overline k}/2}(w_{\overline k}))}
\leq \frac{C}{\rho_{\overline k} \rho_0^n}.
\end{equation}
Now by using interpolation inequality \eqref{inter}, \eqref{4.210} \eqref{1)-10} and \eqref{2)-11},
taking into account that $\overline y\in B_{\rho_{\overline k}/2}(w_{\overline k})$,
we have
\begin{equation}
\label{1-18}
|v(\overline y,\overline t)|\leq
C\left(\frac{\rho_0}{\rho_{\overline k}}\right)^{n/2}H
\left(\left(\frac{\tilde\varepsilon}{H}\right)^{s_5}\right)^{\theta_*^{\overline k}},
\end{equation}
where $\theta_*\in(0,1)$ and $s_5\in(0,1)$ (absolute constant)
and $C$ depends on the a priori data only.
Now evaluating $\overline k$ in terms of $h$ and recalling that
$w(y,s)=\mathcal U(y,s;\overline y,\overline t)=v(\overline y,\overline t)$
we have, for every $(y,s)\in\tilde\Xi_{\rho_0}\times(0,T)$
\begin{equation}
\label{2-33}
\|w\|_{L^\infty(\tilde\Xi_{\rho_0}\times(0,T))}\leq
\frac{C}{h^{n/2}}
\left(\frac{\varepsilon}{C}\right)^{\frac{1}{C}
\left(\frac{h}{\rho_0}\right)^{\frac{|\log\theta_*|}{|\log a|}}},
\end{equation}
where $C$ depends on the a priori data only.
Arguing as above to estimate $|w(y_1,t_1)|$
and recalling that $w(y_1,t_1)=\mathcal U(y_1,t_1;\overline y,\overline t)$
we have
\begin{equation}
\label{1-23}
|\mathcal U(y_1,t_1;\overline y,\overline t)|\leq
\frac{C}{h^n}
\varepsilon^{\frac{1}{C}\left(\frac{h}{\rho_0}\right)^B},
\end{equation}
where $C$ and $B$ depend on the a priori data only.
Finally, using Proposition \ref{prefb}
and proceeding as in Step 1 we obtain \eqref{modulo}.
\cvd

\section{Proof of the Auxiliary Results}

\subsection{Proof of Proposition \ref{dm2}}
\label{distanza}
{\bf Proof of Proposition \ref{dm2}.}
We recall that, for a given a subset $A$ of $\mathbb R^n$, we denote by
$[A]_{\varepsilon}=\left\{x\in\mathbb R^n\,:\,\mathrm{dist}(x,A)\leq\varepsilon\right\},$
$(A)_{\varepsilon}=\left\{x\in A\,:\,\mathrm{dist}(x,\partial A)>\varepsilon\right\}$
and
$[\partial A]_{\varepsilon}=\left\{x\in\mathbb R^n\,:\,\mathrm{dist}(x,\partial A)\leq\varepsilon\right\}$.
We remark that $[A]_\varepsilon\setminus(A)_\varepsilon=[\partial A]_\varepsilon$ and
$(A)_\varepsilon\subset A\subset [A]_\varepsilon$.

\noindent Let $d:=d_{\mathcal H}(\overline D_1,\overline D_2)$ and $r:=d_{\mathcal H}(\partial D_1,\partial D_2)$.
If $d=0$ then \eqref{dm4} holds trivially. Assume $d>0$.
Without loss of generality, we can assume that there exists $\overline x\in\overline D_1$ such that
$d=\mathrm{dist}(\overline x,\overline D_2)$. Since $d>0$ we have that $\overline x\notin\overline D_2$
and therefore $d=\mathrm{dist}(\overline x,\partial D_2)$.
If $\overline x\in\partial D_1$, then \eqref{dm4} is trivially true.
Assume $\overline x\in D_1\setminus\overline D_2$. We have
for every $x\in\overline D_1\setminus\overline D_2$
$$\mathrm{dist}(x,\partial D_2)=\mathrm{dist}(x,\overline D_2)
\leq\mathrm{dist}(\overline x,\overline D_2)=\mathrm{dist}(\overline x,\partial D_2).$$
Thus for every $x\in\overline D_1\setminus\overline D_2$ we have
$$\mathrm{dist}(x,\partial D_2)\leq\mathrm{dist}(\overline x,\partial D_2),$$
that is $\overline x$ is a maximum point in the set $\overline D_1\setminus\overline D_2$
for the function $\mathrm{dist}(\cdot,\partial D_2)$.
In the set $\tilde A=\mathrm{Int}\left([\partial D_2]_{\rho_0/E}\right)\setminus\partial D_2$,
the function $\mathrm{dist}(\cdot,\partial D_2)$ is $C^2$ and
\begin{equation}
\label{3-3}
\left|\nabla_x\mathrm{dist}(x,\partial D_2)\right|>0\qquad \forall\,x\in\tilde A.
\end{equation}
Since $\overline x$ is a maximum point and $\overline x\notin\partial D_2$,
by \eqref{3-3} we have
\begin{equation}
\label{4-3}
\mathrm{dist}(\overline x,\partial D_2)\geq\frac{\rho_0}{E}.
\end{equation}
Otherwise, recalling that $x$ is a maximum point of $\mathrm{dist}(x,\partial D_2)$
interior to $D_1\setminus\overline D_2$, if $\mathrm{dist}(\overline x,\partial D_2)<\rho_0/E$
we should have $\nabla_x \mathrm{dist}(x,\partial D_2)=0$ contradicting \eqref{3-3}.
First let us assume $r$ be such that
\begin{equation}
\label{5-4}
r<\min\left\{\frac{\rho_0}{E},\frac{\rho_0}{2}\right\}.
\end{equation}
We can write $\mathbb R^n=(D_2)_r\cup[\partial D_2]_r\cup(\mathbb R^n\setminus[D_2]_r)$.
By \eqref{4-3} and \eqref{5-4} we have $\overline x\notin[\partial D_2]_r$.
Since $(D_2)_r\subset D_2$ and $\overline x\notin\overline D_2$
we have that $\overline x\in\mathbb R^n\setminus[D_2]_r$.
Recalling that $r<\rho_0/2$ and $\mathbb R^n\setminus D_2$ is connected,
we have that $\mathbb R^n\setminus[D_2]_r$ is connected. Thus there exists a continuous path
$$\gamma\,:\,[0,1)\to\mathbb R^n$$
such that
\begin{subequations}
\label{6-5}
\begin{equation}
\label{6-1-5}
\gamma([0,1))\subset\mathbb R^n\setminus[D_2]_r,
\end{equation}
\begin{equation}
\label{6-2-5}
\gamma(0)=\overline x\qquad \lim\limits_{t\to1^-}\gamma(t)=\infty.
\end{equation}
\end{subequations}
Since $d_{\mathcal H}(\partial D_1,\partial D_2)=r$ and
$\partial D_1\subset[\partial D_2]_r\subset[D_2]_r$, by \eqref{6-1-5} we have
$$\gamma([0,1))\cap\partial D_1=\emptyset$$
which is a contradiction since $\overline x\in D_1$ and $D_1$ is bounded.
Thus we cannot connect $\overline x$ and $\infty$ with a path that does not
intersect $\partial D_1$. Hence $\overline x\in\partial D_1$.
Thus if
$d_{\mathcal H}(\partial D_1,\partial D_2)\leq\delta\rho_0$,
with $\delta=\min\left\{\frac{1}{E},\frac{1}{2}\right\}$,
\eqref{5-4} is satisfied and we have
\begin{equation*}
d_{\mathcal H}(\overline D_1,\overline D_2)=\mathrm{dist}(x,\partial D_2)
\leq r.
\end{equation*}
On the other side, if $d_{\mathcal H}(\partial D_1,\partial D_2)>\delta\rho_0$
we have trivially
\begin{equation}
\label{3-1U}
d_{\mathcal H}(\overline D_1,\overline D_2)\leq 2\textrm{diam}(\Omega)
\leq\frac{2\textrm{diam}(\Omega)}{\delta\rho_0} d_{\mathcal H}(\partial D_1,\partial D_2)
\end{equation}
and the proposition is proven.
\cvd

\subsection{Proof of Proposition \ref{5-16Tpr} and Asymptotic Estimates for the Fundamental Solution}
\label{fax-sec}

We shall make use of the following regularity theorem,
whose proof can be found in \cite{la2}, \cite[Ch. III, Sec. 13]{la}.

\begin{teo}
\label{L}
Let $\lambda$, $M$ and $r$ be positive numbers with $\lambda\in(0,1]$.
Let $u\in H^{1,\frac{1}{2}}(B'_r\times(-r,r)\times(-r^2,r^2))$
be solution to
\begin{equation}
\label{1-1}
\mathrm{div}\left(A(x,t)\nabla_xu\right)+b(x,t)\cdot\nabla_xu-
\partial_tu=0,
\end{equation}
where $A(x,t)$ and $b(x,t)$ are respectively a symmetric $n\times n$
matrix and a vector valued function satisfying the following conditions
\begin{subequations}
\begin{equation}
\label{2-1}
\lambda|\xi|^2\leq A(x,t)\xi\cdot\xi
\leq\lambda^{-1}|\xi|^2,
\end{equation}
for all $(x,t)\in B'_r\times(-r,r)\times(-r^2,r^2)$ and for all $\xi\in\mathbb R^n$,
\begin{multline}
\label{3-1}
r\sum_{i=1}^n\|\partial_i A\|_{L^\infty(B'_r\times(-r,0)\times(-r^2,r^2))}
+r\sum_{i=1}^n\|\partial_i A\|_{L^\infty(B'_r\times(0,r)\times(-r^2,r^2))}\\
r^2\|\partial_tA\|_{L^\infty(B'_r\times(-r,r)\times(-r^2,r^2))}\leq M,
\end{multline}
\begin{equation}
\label{4-1}
r\|b\|_{L^\infty(B'_r\times(-r,r)\times(-r^2,r^2))}\leq M.
\end{equation}
\end{subequations}
Then there exist positive constants $\beta\in(0,1)$ and
$C$ such that for every $\rho<\frac{r}{2}$ and all
$(x,t)\in B'_{r-2\rho}\times(-(r-2\rho),(r-2\rho))
\times(-r^2+4\rho^2,r^2)$ the following inequality holds
\begin{eqnarray}
\label{1-2}
&&\rho\|\nabla_xu\|_{L^\infty(B'_\rho(x')\times(-\rho+x_n,\rho+x_n)
\times(-\rho^2+t,t))}\\[2mm]
&&+\rho^{\beta+1}[\nabla_x u]_{\beta;(B'_\rho(x')\times(-\rho+x_n,\rho+x_n)
\times(-\rho^2+t,t))
\cap (B'_\rho\times(-r,0)\times(-r^2,r^2))}\nonumber\\[2mm]
&&+\rho^{\beta+1}[\nabla_x u]_{\beta;(B'_\rho(x')\times(-\rho+x_n,\rho+x_n)
\times(-\rho^2+t,t))
\cap (B'_\rho\times(0,r)\times(-r^2,r^2))}\nonumber\\[2mm]
&\leq&\frac{C}{\rho^{\frac{n}{2}+1}}\left\{
\int_{B'_{2\rho}(x')\times(-2\rho+x_n,2\rho+x_n)\times(-4\rho^2+t,t)}
u^2(\xi,\tau)d\xi d\tau\right\}^{1/2}.\nonumber
\end{eqnarray}
Here $\beta$ depends on $n$ only and $C$ depends on
$\lambda$, $M$ and $n$ only.
\end{teo}

Before proving Proposition \ref{5-16Tpr} we give
the following estimate which is needed in the proof.
We recall  that $Q_\rho(x_0,t_0)=B_r(x_0)\times(t_0-\rho^2,t_0)$.

\begin{prop}
\label{fax}
There exist constant $C\geq1$ and $0<\delta_1<1$
depending on $k$ and $n$ only such that
the following inequality holds.
\begin{equation}
\int_{Q_\rho(x_0,t_0)}|\Gamma(x,t;\xi,\tau)|^2dx\,dt
\leq
C\frac{\rho^n}{(t_0-\tau)^{n-1}}
\mathrm e^{-\frac{|x_0-\xi|^2}{C(t_0-\tau)}},
\end{equation}
where $\rho=\delta_1[|x_0-\xi|^2+t_0-\tau]^{1/2}$.
\end{prop}
\textbf{Proof.}
From the inequality \eqref{av3.3} we have
\begin{equation}
\label{fax1-2}
\int_{Q_\rho(x_0,t_0)}|\Gamma(x,t;\xi,\tau)|^2dx\,dt
\leq C\int_{Q_\rho(x_0,t_0)}
\frac{1}{(t-\tau)^{n}}
\mathrm e^{-\frac{|x-\xi|^2}{C_1(t-\tau)}}
\chi_{[\tau,+\infty)}dx\,dt,
\end{equation}
where $C_1$ depends on $k$ and $n$ only.
In what follows we denote by $I$ the integral
at the right-hand side of \eqref{fax1-2}.
We distinguish two cases

i) $t_0-\rho^2<\tau<t_0$,

ii) $\tau<t_0-\rho^2$.

\noindent Let us consider case i). It is easy to see that
there exists an absolute constant $C\geq1$ such that
\begin{equation}
\label{fax2-2}
C^{-1}\rho\leq|x-\xi|\leq C\rho\qquad\forall\,x\in B_\rho(x_0).
\end{equation}
By \eqref{fax2-2} we have
\begin{equation}
\label{fax3-2}
I\leq c_n\rho^n\int_0^{t_0-\tau}
s^{-n}\mathrm e^{-\frac{\rho^2}{C_2s}}ds,
\end{equation}
where $c_n$ is an absolute constant depending on $n$
only and $C_2$ depends on $k$ and $n$ only.
Now if $0<t_0-\tau<\frac{\rho^2}{nC_2}$,
being $s\to s^{-n}\mathrm e^{-\frac{\rho^2}{C_2s}}$
an increasing function in $(0,\frac{\rho^2}{nC_2})$,
by \eqref{fax3-2} we get
\begin{equation}
\label{fax4-2}
I\leq\frac{\rho^n}{(t_0-\tau)^{n-1}}
\mathrm e^{-\frac{\rho^2}{C(t_0-\tau)}}.
\end{equation}
Otherwise, if $\frac{\rho^2}{nC_2}<t_0-\tau<\rho^2$
then
since
$$\max _{(0,+\infty)}\{s^{-n}\mathrm e^{-\rho^2/(C_2s)}\}
=\frac{(nC_2)^n}{\rho^{2n}}\mathrm e^{-1/n}$$
and now $t-\tau$ is of the same order of $\rho^2$ we have
\begin{eqnarray*}
&&\rho^{-n}I\leq
C_1\int_0^{\rho^2}s^{-n}\mathrm e^{-\frac{\rho^2}{C_2s}}ds
\leq \frac{C}{(t_0-\tau)^{n-1}}
\mathrm e^{-\frac{\rho^2}{C_2(t_0-\tau)}}.
\end{eqnarray*}
By the last inequality and \eqref{fax4-2} we get the Proposition
in case i).

Let us consider now case ii).
It is easy to see that
\begin{equation}
\label{fax1-3}
6\rho^2\leq|x-\xi|^2+t-\tau\leq60\rho^2,
\end{equation}
for every $(x,t)\in Q_\rho(x_0,t_0)$.
Moreover, denoting
$$M_\rho=\max\left\{\frac{\mathrm e^{-\frac{|x-\xi|^2}{C_1(t-\tau)}}}
{(t-\tau)^n}\,:\,(x,t)\in Q_\rho(x_0,t_0)\right\}$$
and taking into account \eqref{fax1-3} we get
\begin{equation}
\label{fax2-3}
M_\rho\leq C\left(\frac{C_1}{\rho^2}\right)^n,
\end{equation}
where $C$ depends on $n$ only.
Now, since $\tau< t_0-\rho^2$ we have
\begin{equation}
\label{fax1-4}
\frac{|x_0-\xi|^2}{t_0-\tau}\leq4.
\end{equation}
Therefore by \eqref{fax2-3} and \eqref{fax1-4} we
get the Proposition in case ii) as well.
\cvd

\bigskip

{\bf Proof of Proposition \ref{5-16Tpr}.}
Let $U$ be a solution of the equation $\mathcal LU=0$,
where $\mathcal L=\partial_t-\mathrm{div}(1+(k-1)\chi_Q \nabla)$.
We recall the following regularity estimate (see \cite{la})
\begin{equation}
\label{3-5}
\|\nabla U\|_{L^\infty(Q^{\pm}_r(\overline x,\overline t))}
\leq\frac{C}{r^{\frac{n+4}{2}}}
\left(\int_{Q_{2r}(\overline x,\overline t)}U^2(x,t)dxdt\right)^{1/2},
\end{equation}
where $Q^+_r(\overline x,\overline t)=Q_r(\overline x,\overline t)\cap Q$
(we recall $Q=D(\mathbb R)$)
and $Q^-_r(\overline x,\overline t)=Q_r(\overline x,\overline t)\setminus Q^+(\overline x,\overline t)$.
Applying \eqref{3-5} to the function $\Gamma(\cdot,\cdot;\xi,\tau)$ we get
\begin{equation}
\label{2-6}
\|\nabla \Gamma(\cdot,\cdot;\xi,\tau)\|_{L^\infty(Q_{\rho}^{\pm}(x_0,t_0))}
\leq\frac{C}{\rho^{\frac{n+4}{2}}}
\left[\int_{Q_{2\rho}(x_0,t_0)}|\Gamma(x,t;\xi,\tau)|^2dxdt\right]^{1/2},
\end{equation}
where
\begin{equation}
\label{1-6}
\rho=\frac{1}{4}\left[|x_0-\xi|^2+t_0-\tau\right]^{1/2}.
\end{equation}
Applying Proposition \ref{fax} to the right hand side of \eqref{2-6} we have
\begin{equation*}
\|\nabla\Gamma(\cdot,\cdot;\xi,\tau)\|_{L^\infty(Q_{\rho}^{\pm}(x_0,t_0))}
\leq\frac{C}{\rho^{\frac{n+4}{2}}}
\left[\frac{\rho^{n}}{(t_0-\tau)^{n-1}}
\,\mathrm e^{-\frac{|x_0-\xi|^2}{C(t_0-\tau)}}\right]^{1/2}.
\end{equation*}
Since
$$\frac{1}{16\rho^2}\leq\frac{1}{t_0-\tau}$$
we obtain \eqref{5-16T}.
\cvd

\bigskip

In order to state the next theorem we introduce some notations.
Let
$\varphi:B_{\rho_0}'\times(-\rho_0^2,\rho_0^2)\to\mathbb R$ such
that it is differentiable with respect to $t$ and $x_i$,
$i=1,\dots,n-1$, it is twice differentiable with respect to $x_i$,
$i=1,\dots,n-1$, and $\partial_t\varphi$ is differentiable with
respect to  $x_i$, $i=1,\dots,n-1$. We assume that
\begin{equation}
\label{1-1-12}
\varphi(0,0)=|\nabla_{x'}\varphi(0,0)|=0
\end{equation}
and
\begin{multline}
\label{2-1-12}
\rho_0^2\|D_{x'}^2\varphi\|_{L^\infty(B'_{\rho_0}\times(-\rho^2,\rho^2))}
+\rho_0^2\|\partial_t\varphi\|_{L^\infty(B'_{\rho_0}\times(-\rho^2,\rho^2))}
\\[2mm]
+\rho_0^3\|\partial_t\nabla_{x'}
\varphi\|_{L^\infty(B'_{\rho_0}\times(-\rho^2,\rho^2)}
\leq E\rho_0.
\end{multline}
We shall denote by
$$\tilde Q_{\varphi,\rho_0}^+=\{x\in B_{\rho_0}
\times(-\rho_0^2,\rho_0^2)\,:\,x_n>\varphi(x',t)\},$$
and by $\Gamma_{\tilde Q_{\varphi,\rho_0}^+}(x,t;y,s)$ the fundamental solution of
the operator $\partial_t-\mathrm{div}((1+(k-1)\chi_{\tilde Q_{\varphi,\rho_0}^+}\nabla)$,
that is
\begin{multline*}
\partial_t \Gamma_{\tilde Q_{\varphi,\rho_0}^+}(x,t;y,s)\\
-\mathrm{div}\left((1+(k-1)\chi_{\tilde Q_{\varphi,\rho_0}^+})
\nabla\Gamma_{\tilde Q_{\varphi,\rho_0}^+}(x,t;y,s)\right)
=-\delta(x-y,t-s),
\end{multline*}
where $(y,s)\in\mathbb R^{n+1}$.

\begin{teo}[Asymptotic Estimate]
\label{STIMA ASINTOTICA}
Let $\varphi$ and $\Gamma_{\tilde Q_{\varphi,\rho_0}^+}(x,t;y,s)$ as above.
Then there exists a constant $C\geq1$ depending on $n$ and $E$
only such that
\begin{eqnarray}
\label{STIMA ASINTOTICA1}
&&\left|\Gamma_{\tilde Q_{\varphi,\rho_0}^+}(x,t;y,0)
-\Gamma_+(x,t;y,0)\right|
\leq C\frac{[|x-y|^2+t]^{1/2}}{\rho_0}
\frac{\mathrm e^{-\frac{|x-y|^2}{Ct}}}{t^{n/2}},\\[3mm]
\label{STIMA ASINTOTICA2}
&&\left|\nabla_x\Gamma_{\tilde Q_{\varphi,\rho_0}^+}(x,t;y,0)
-\nabla_x\Gamma_+(x,t;y,0)\right|\\[2mm]
&&\qquad\qquad\leq C\frac{[|x-y|^2+t]^{\frac{1}{2}(-1+\frac{\beta}{\beta+1})}}
{\rho_0^{\frac{\beta}{1+\beta}}}
\frac{\mathrm e^{-\frac{|x-y|^2}{Ct}}}{t^{n/2}},\nonumber
\end{eqnarray}
where $\beta$ is the one defined in Theorem~\ref{L},
depending on $n$ only, for all
\begin{equation*}
(x,t)\in
\tilde Q_{\varphi,\frac{\rho_0}{C}}
\cap\left\{(x,t)\in\mathbb R^{n+1}\,:\,t>0,\,x_n>\frac{1}{C\rho_0}(|x'|^2+t)\right\}
\end{equation*}
and $y=y_ne_n$, $y_n\in(-\rho_0/C,0)$.
\end{teo}

\begin{oss}
Theorem \ref{STIMA ASINTOTICA} provides an asymptotic estimate for
the fundamental solution $\Gamma(x,t;y,s)$ when $(x,t)$ and $(y,s)$
stay on opposite sides of the interface (given by the graphic $x_n=\varphi(x',t)$).
Our crucial requirement is
that $(y,s)$ approaches the interface in a nontangential way.
\end{oss}

{\bf Proof. of Theorem \ref{STIMA ASINTOTICA}.}
Let $\theta$ be a $C^{\infty}$ function on $\mathbb R$ such that
$0\leq\theta\leq1$, $\theta(s)=0$, for every $s\in\mathbb R\setminus(-2,2)$,
$\theta(s)=1$ for every $s\in(-1,1)$ and $|\theta'(s)|\leq2$ for every
$s\in\mathbb R$.

We define new variables by  $(\xi,\tau)=\Psi(x,t)$,
where $\Psi(x,t)=(\Phi(x,t),t)$ and
$$\left\{\begin{array}{l}
\xi'=x',\\
\xi_n=x_n-\varphi(x',t)\theta\left(\frac{|x'|}{r_1}\right)\theta\left(\frac{x_n}{r_1}\right)
\theta\left(\frac{t}{r^2_1}\right),\\
\tau=t,
\end{array}\right.$$
where $r_1=\rho_0\min\{\frac{1}{4},\frac{1}{32E}\}$.

Sometimes, for the sake of brevity for a fixed $t\in(-\rho_0^2,\rho_0^2)$
we denote by $\Phi^{(t)}(\cdot)$ the map $\Phi(\cdot,t)$
and by $G^{(t)}$ the graph of $\varphi(\cdot,t)$.
It is not difficult to check that $\Psi$ and $\Phi$
have the same regularity properties of $\varphi$ and they are
diffeomorphisms (that preserve orientation) of $\mathbb R^{n+1}$
and $\mathbb R^n$ respectively.
We denote by $\Phi^{-1}(\cdot,t)$ the inverse of $\Phi^{(t)}(\cdot)$.
The following properties hold:
\begin{subequations}
\label{phi}
\begin{eqnarray}
\label{a}
&&\Phi^{(t)}\left(G^{(t)}\cap(B'_{r_1}\times(-r_1,r_1)) \right)=
\{x\in B'_{r_1}\times(-r_1,r_1)\,:\,x_n=0\},\\[2mm]
\label{b}
&&\Psi(x,t)=(x,t),\\ [1mm]\nonumber
&&\qquad\forall\,
(x,t)\in\mathbb R^{n+1}\setminus \left(\left(B'_{2r_1}\times(-2r_1,2r_1)\right)
\times(-2r_1^2,2r_1^2)\right),
\\[2mm]
\label{c}
&&C^{-1}|x_1-x_2|\leq|\Phi^{(t)}(x_1)-\Phi^{(t)}(x_2)|\leq C|x_1-x_2|,\quad \forall\,
x_1,x_2\in\mathbb R^n,\\[2mm]
\label{d}
&&|\Phi^{(t)}(x)-x|\leq \frac{C}{\rho_0}|x|^{2},\quad
\forall\,x\in\mathbb R^n,\\[2mm]
\label{e}
&&|D_x\Phi^{(t)}(x)-I|\leq\frac{C}{\rho_0}|x|,\quad\forall\,
x\in\mathbb R^n,
\end{eqnarray}
\end{subequations}
where $C$, $C\geq1$, depends on $E$ only,
$I$ denotes the identity matrix and
$D_x\Phi^{(t)}$ is the jacobian matrix with respect to variable $x$.
For $y_n\in(-\frac{r_1}{2},0)$ and $\sigma\in(-r_1^2,r_1^2)$
we denote $y=y_n e_n$ and $\eta=\Phi(y,\sigma)$.
Furthermore we shall use the following notation
$\tilde{\Gamma}(\xi,\tau;\eta,\sigma)
=\Gamma(\Psi^{-1}(\xi,\tau);\Psi^{-1}(\eta,\sigma))$,
and $\gamma(\xi,\tau)=\det J(\xi,\tau)$,
where $J(\xi,\tau)=(D_x\Phi)(\Psi^{-1}(\xi,\tau))$.
We have that $\tilde\Gamma(\xi,\tau;\eta,\sigma)$
is a solution to
\begin{equation}
\label{BC}
\mathrm{div}\left(\tilde B(\xi,\tau)\nabla_{\xi}\tilde\Gamma\right)
+C(\xi,\tau)\nabla_{\xi}\tilde\Gamma-\partial_{\tau}\tilde\Gamma
=-\gamma(\eta,\sigma)\delta(\xi-\eta,\tau-\sigma),
\end{equation}
where $\tilde B(\xi,\tau)=(1+(k-1)\chi^+)B(\xi,\tau)$,
$B(\xi,\tau)=\left(J(\xi,\tau)\right)\left(J(\xi,\tau)\right)^*$
and
$C(\xi,\tau)=J(\xi,\tau)\frac{\partial\Phi^{-1}(\xi,\tau)}{\partial\tau}
-\frac{\tilde B(\xi,\tau)}{\gamma(\xi,\tau)}\nabla_{\xi}\gamma(\xi,\tau)$.

Since we want to study the asymptotic behaviour of $\tilde\Gamma(\xi,\tau;\eta,0)$,
we shall denote $\tilde\Gamma(\xi,\tau;\eta,0)$ by $\tilde\Gamma(\xi,\tau;\eta)$.

By \eqref{phi}, we have that
\begin{equation}
\label{1N-3}
B(0,0)=I\qquad\textrm{and}\qquad
\Vert B\Vert_{L^{\infty}(\Omega\times(0,T))}
+\rho_0[B]_{1,\Omega\times(0,T)}\leq C,
\end{equation}
where $C$ depends on $E$ only. Denote by
\begin{equation}
\label{defR}
R(\xi,\tau;\eta)=\tilde\Gamma(\xi,\tau;\eta)-\gamma(\eta,0)\Gamma_+(\xi,\tau;\eta),
\end{equation}
where $\Gamma_+(\xi,\tau;\eta)=\Gamma_+(\xi,\tau;\eta,0)$ is the fundamental solution
to the operator $\mathrm{div}\left((1+(k-1)\chi^+)\nabla_{\xi}\right)
-\partial_\tau$. We have
$$\mathrm{div}\left((1+(k-1)\chi^+)\nabla_{\xi}R\right)-\partial_{\tau}R=
F(\xi,\tau;\eta),$$
where
\begin{eqnarray*}
F(\xi,\tau;\eta)&=&-C(\xi,\tau)\nabla_{\xi}\tilde\Gamma(\xi,\tau;\eta)\\[2mm]
&&+\mathrm{div}\left((1+(k-1)\chi^+)(I-B(\xi,\tau))\nabla_{\xi}\tilde\Gamma(\xi,\tau;\eta)\right),
\end{eqnarray*}
notice that $F(\xi,\tau;\eta)=0$ if $\tau<0$,
\begin{equation}
\label{1N-4}
R(\xi,\tau;\eta)=0,\quad\textrm{for } \tau<0.
\end{equation}
Therefore, \cite{ar}
$$R(\xi,\tau;\eta)=\int_0^{\tau}\int_{B_{2r_1}}F(\zeta,s;\eta)
\Gamma_+(\xi,\tau;\zeta,s)d\zeta ds,\qquad\textrm{if }\tau>0.$$
We have
\begin{equation}
\label{stimaR}
|R(\xi,\tau;\eta)|\leq J_1+J_2,
\end{equation}
where
\begin{subequations}
\begin{equation}
\label{J1}
J_1=\left|\int_0^{\tau}\int_{B_{2r_1}}C(\zeta,s)
\nabla_{\zeta}\tilde\Gamma(\zeta,s;\eta)\Gamma_+(\xi,\tau;\zeta,s)d\zeta ds\right|
\end{equation}
and
\begin{multline}
\label{J2}
J_2=\left|\int_0^{\tau}\int_{B_{2r_1}}\left(1+(k-1)\chi^+\right)
\left(I-B(\zeta,s)\right)\right.\\
\times\left.\nabla_{\zeta}\tilde\Gamma(\zeta,s;\eta)\cdot
\nabla_\zeta\Gamma_+(\xi,\tau;\zeta,s)d\zeta ds\right|.
\end{multline}
\end{subequations}
By Proposition \ref{5-16Tpr} and Lemma \ref{fr} we have
\begin{equation}
\label{J1stima}
J_1\leq\frac{C}{\rho_0\tau^{\frac{n-1}{2}}}\textrm e^{-\frac{|\xi-\eta|^2}{C_1\tau}}.
\end{equation}
where $C$, $C_1$, $C\geq1$, $C_1\geq1$, depend on $E$ only.
By Proposition \ref{5-16Tpr} and \eqref{1N-3} we get
\begin{equation}
J_2\leq J_{2,1}+J_{2,2},
\end{equation}
with
\begin{equation}
J_{2,1}=\frac{C}{\rho_0}\int_0^{\tau}\int_{\mathbb R^n}
s^{-\frac{n}{2}}\textrm e^{-\frac{|\zeta-\eta|^2}{C_1 s}}
(\tau-s)^{-\frac{n+1}{2}}
\textrm e^{-\frac{|\xi-\zeta|^2}{C_1(\tau-s)}}d\zeta ds,
\end{equation}
and
\begin{equation}
J_{2,2}=\frac{C}{\rho_0}\int_0^{\tau}\int_{\mathbb R^n}
|\zeta|s^{-\frac{n+1}{2}}\textrm e^{-\frac{|\zeta-\eta|^2}{C_1 s}}
(\tau-s)^{-\frac{n+1}{2}}
\textrm e^{-\frac{|\xi-\zeta|^2}{C_1(\tau-s)}}d\zeta ds,
\end{equation}
where $C,C_1$, $C\geq1$, $C_1\geq 1$, depend on $E$ only.

By Lemma~\ref{fr} we obtain
\begin{equation}
\label{J21stima}
J_{2,1}\leq C\frac{\tau^{-\frac{n}{2}+\frac{1}{2}}}{\rho_0}
\textrm e^{-\frac{|\xi-\eta|^2}{C_1\tau}}.
\end{equation}
Let us consider now $J_{2,2}$.
Performing a change of variables we get
\begin{eqnarray*}
J_{2,2}&=&\frac{C}{\rho_0}\left(\frac{\sqrt{C_1}}{2}\right)^n
\frac{\textrm e^{-\frac{|\xi-\eta|^2}{C_1\tau}}}{\tau^{n/2}}\\[2mm]
&&\qquad\times\int_{0}^{1}\int_{\mathbb R^n}\left|
\frac{C_1}{2}(\tau(1-\lambda)\lambda)^{1/2}z+\lambda(\xi-\eta)+\xi\right|
\frac{\textrm e^{-|z|^2}}{\sqrt{(1-\lambda)\lambda}}dzd\lambda\\[2mm]
&\leq&\frac{C}{\rho_0}
\frac{\mathrm e^{-\frac{|\xi-\eta|^2}{C_1\tau}}}{\tau^{n/2}}
\left[(|\xi|^2+\tau)^{1/2}+|\xi-\eta|\right],
\end{eqnarray*}
where $C$ depends on $E$ only.
Now, denoting by $C_2=\max\limits_{s\in(0,+\infty)}s^{1/2}\mathrm e^{-\frac{1}{2C_1s}}$,
we have
\begin{eqnarray*}
&&\frac{\mathrm e^{-\frac{|\xi-\eta|^2}{C_1\tau}}}{\tau^{n/2}}|\xi-\eta|
=\frac{1}{\tau^{\frac{n-1}{2}}}
\left(\frac{|\xi-\eta|^2}{\tau}\right)^{1/2}\mathrm e^{-\frac{|\xi-\eta|^2}{2C_1\tau}}
\mathrm e^{-\frac{|\xi-\eta|^2}{2C_1\tau}}
\leq\frac{C_2}{\tau^{\frac{n-1}{2}}}
\mathrm e^{-\frac{|\xi-\eta|^2}{2C_1\tau}}.
\end{eqnarray*}
Thus
\begin{equation}
\label{1N-7}
J_{2,2}
\leq\frac{C}{\rho_0}
\frac{\mathrm e^{-\frac{|\xi-\eta|^2}{2C_1\tau}}}{\tau^{n/2}}(|\xi|^2+\tau)^{1/2}.
\end{equation}
Now since $\eta=e_n\eta_n$, $\eta_n<0$ and $\xi_n>0$, we have
$|\xi-\eta|^2=|\xi|^2-2\eta_n\xi_n+|\eta|^2\geq|\xi|^2$.
Such an inequality and \eqref{1N-4}, \eqref{J1stima}, \eqref{J21stima}, \eqref{1N-7}
give
\begin{equation}
\label{1-4}
|R(\xi,\tau;\eta)|\leq\frac{C}{\rho_0}
\chi_{_{\mathbb R^n\times[0,+\infty)}}
\frac{\mathrm e^{-\frac{|\xi-\eta|^2}{C\tau}}}{\tau^{n/2}}
(|\xi-\eta|^2+\tau)^{1/2},
\end{equation}
for every $\xi\in B^+_{2r_1}$ and $\tau\in(0,4r_1^2)$,
where $C$, $C\geq1$, depends on $E$ only.
Let $\delta_1$ be the constant defined in Proposition \ref{fax}
($\delta_1\in(0,1)$) and, for fixed $\overline \xi\in B_{r_1/8}^+$,
$\eta_n\in(-r_1/8,0)$, $\eta=e_n\eta_n$, $\overline\tau\in(0,(r_1/8)^2)$
denote by
$$h=\frac{\delta_1}{4}[|\overline\xi-\eta|^2+\overline\tau]^{1/2}.$$
We have
$$\mathrm{div}\left(\tilde B(\xi,\tau)\nabla_\xi\tilde\Gamma\right)
+C(\xi,\tau)\nabla_\xi\tilde\Gamma-\partial_\tau\tilde\Gamma=0,$$
in $B'_{h/2}(\overline\xi')\times
(\overline\xi_n-h/2,\overline\xi_n+h/2)
\times(\overline\tau-(h/2)^2,\overline\tau]$, where
$\tilde B$ and $C$ are defined above. Therefore
by Theorem \ref{L} and Proposition \ref{fax} we get
\begin{equation}
\label{1-9}
[\nabla_\xi\tilde\Gamma(\cdot,\overline\tau,\eta)]_{\beta,Q}
\leq \frac{C}{h^{2+\beta}}\frac{1}{\overline\tau^{\frac{n-1}{2}}}
\mathrm e^{-\frac{|\overline\xi-\eta|^2}{C\overline\tau}},
\end{equation}
where where $Q=B'_{h/4}(\overline\xi')\times(\overline\xi_n,\overline\xi_n+h/4)$.
Since a similar inequality holds true for $\nabla_\xi\Gamma_+(\cdot,\overline\tau;\eta)$,
by \eqref{defR} we obtain
\begin{equation}
\label{1-10}
[\nabla_\xi R(\cdot,\overline\tau;\eta)]_{\beta;Q}
\leq \frac{C}{h^{2+\beta}}\frac{1}{\overline\tau^{\frac{n-1}{2}}}
\mathrm e^{-\frac{|\overline\xi-\eta|^2}{C\overline\tau}}.
\end{equation}
In \eqref{1-9} and \eqref{1-10}, $C$, $C\geq1$, depends on $E$ only.
Now we recall the following interpolation inequality
\begin{equation}
\label{2-10}
\|\nabla f\|_{L^\infty(Q)}\leq
C\left(\|f\|_{L^\infty(Q)}^{\frac{\beta}{1+\beta}}
|\nabla f|_{\beta;Q}^{\frac{1}{1+\beta}}+
\frac{1}{h}\|f \|_{L^\infty(Q)}\right).
\end{equation}
Since \eqref{1-4} easily yields
\begin{equation}
\label{1star}
\|R(\cdot,\overline\tau;\eta)\|_{L^\infty(Q)}\leq
\frac{C}{\rho_0}\frac{\mathrm e^{-\frac{|\overline\xi-\eta|^2}{C\overline\tau}}}
{\overline\tau^{n/2}}h,
\end{equation}
where $C$, $C\geq1$, depends on $E$ only, we obtain
by \eqref{1-10} and \eqref{2-10}
\begin{equation}
\label{1-c1}
|\nabla_\xi R(\overline\xi,\overline\tau;\eta)|\leq
\frac{C}{\rho_0}\left(\frac{h}{\rho_0}\right)^{-1+\frac{\beta}{\beta+1}}
\frac{\mathrm e^{-\frac{|\overline\xi-\eta|^2}{C_1\overline\tau}}}
{\overline\tau^{n/2}},
\end{equation}
for every $\overline\xi\in B^+_{r_1/8}$, $\eta=e_n\eta_n$,
$\eta_n\in (-r_1/8,0)$, $\overline\tau\in(0,(r_1/8)^2]$,
where $C$, $C\geq1$, depends on $E$ only.

Let us go back to the original coordinates $(x,t)$.
First of all let us estimate the function $g$ defined by
\begin{equation}
\label{1N-9}
g(x,t;y):=R(\Phi^{(t)}(x),t;\Phi^{(0)}(y))=
R(\Phi^{(t)}(x),t;e_ny_n).
\end{equation}
To carry out the estimates, up to the end of the proof, we always consider
$x$ and $y_n$ such that $x\in B_{\delta\rho_0}(\delta\rho_0e_n)$,
$y_n\in(-\delta\rho_0,0)$, where $\delta$, $\delta\in(0,1)$, may change
from line to line, but it shall depend on $E$ only.
Notice that for every $x\in B_{\delta\rho_0}(\delta\rho_0e_n)$
we have $x_n>0$. Also notice that
\begin{equation}
\label{N1-10}
|x|\leq|x-e_ny_n|,\quad x\in B_{\delta\rho_0}(\delta\rho_0e_n),\,y_n\in(-\delta\rho_0,0).
\end{equation}
By such an inequality and  \eqref{d} we have
\begin{equation}
\label{star2}
|\Phi^{(t)}(x)-x|\leq\frac{C}{\rho_0}|x-e_ny_n|^2
\end{equation}
for $x\in B_{\delta\rho_0}(\delta\rho_0e_n)$, $y_n\in(-\delta\rho_0,0)$
where $C$ depends on $E$ only.
By \eqref{star2}, \eqref{N1-10} and the triangle inequality we have
\begin{equation}
\label{2-c3}
C^{-1}|x-e_ny_n|\leq|\Phi^{(t)}(x)-e_ny_n|\leq C|x-e_ny_n|,
\end{equation}
for $x\in B_{\delta\rho_0}(\delta\rho_0e_n)$, $y_n\in(-\delta\rho_0,0)$
where $C$, $C\geq1$, depends on $E$ only.
By \eqref{1star}, \eqref{1-c1}, \eqref{e}, \eqref{1N-9},
\eqref{N1-10}, \eqref{2-c3} we obtain
\begin{equation}
\label{1-41}
|g(x,t,y)|\leq C\frac{\mathrm e^{-\frac{|x-e_ny_n|^2}{Ct}}}{t^{n/2}}
\left[\frac{|x-e_ny_n|^2+t}{\rho_0^2}\right]^{1/2},
\end{equation}
and
\begin{equation}
\label{2-c5}
|\nabla_x g(x,t;y)|\leq
\frac{C}{\rho_0}\frac{\mathrm e^{-\frac{|x-e_ny_n|^2}{Ct}}}
{t^{n/2}[|x-e_ny_n|^2+t]^{1/2}}
\left(\frac{[|x-e_ny_n|^2+t]^{1/2}}{\rho_0}
\right)^{-1+\frac{\beta}{1+\beta}},
\end{equation}
for every $x\in B_{\delta\rho_0}(\delta\rho_0e_n)$, $y_n\in(-\delta\rho_0,0)$
where $C$, $C\geq1$, depends on $E$ only.
Recalling the definition of $g$ we have that
\begin{multline}
\label{1-c5}
\gamma\left(\Phi^{(0)}(y)\right)
\left(\Gamma(x,t;y,0)-\Gamma_+(x,t;y,0)\right)
=\\g(x,t;y)-\left(1-\gamma(\Phi^{(0)}(y))\right)\Gamma_+(x,t;y,0)\\
-\gamma\left(\Phi^{(0)}(y)\right)\left(\Gamma_+(x,t;y,0)-\Gamma_+(\Phi^{(t)}(x),t;y,0)\right).
\end{multline}
Now for $x\in B_{\delta\rho_0}(\delta\rho_0e_n)$, $y_n\in(-\delta\rho_0,0)$
we have
\begin{equation}
\label{1N-11}
|y|\leq|x-e_ny_n|,
\end{equation}
so such an inequality, \eqref{1N-3}, \eqref{c} and \eqref{av3.3} give
\begin{equation}
\label{2N-11}
\left|\left(1-\gamma(\Phi^{(0)}(y))\right)\Gamma_+(x,t;y,0)\right|
\leq C\left(\frac{|x-e_ny_n|^2+t}{\rho_0^2}\right)^{1/2}
\frac{\mathrm e^{-\frac{|x-e_ny_n|^2}{Ct}}}{t^{n/2}},
\end{equation}
for every $x\in B_{\delta\rho_0}(\delta\rho_0e_n)$, $y_n\in(-\delta\rho_0,0)$
where $C$, $C\geq1$, depends on $E$ only.
In order to estimate from above the third term at the right hand side of \eqref{1-c5}
we use the mean value theorem. By such a theorem, \eqref{1N-3} and Proposition \ref{5-16Tpr}
we get
\begin{multline}
\label{1N-12}
\left|\gamma(\Phi^{(0)}(y))\left(\Gamma_+(\Phi^{(t)}(x),t;y,0)\right)
-\Gamma_+(x,t;y,0)\right|\\
\leq C|x-\Phi^{(t)}(x)|\frac{\mathrm e^{-\frac{|\overline x-e_ny_n|^2}{Ct}}}{t^{\frac{n+1}{2}}},
\end{multline}
where $\overline x=x+\lambda(\Phi^{(t)}(x)-x)$ for a suitable $\lambda\in(0,1)$
and $C$, $C\geq1$, depends on $E$ only.
Now, by triangle inequality, \eqref{d}, \eqref{N1-10} we have
\begin{multline}
\label{1N-13}
|\overline x-e_ny_n|\geq|x-e_ny_n|-|x-\overline x|\\
\geq|x-e_ny_n|\left(1-\frac{C}{\rho_0}|x|\right)\geq
\frac{1}{2}|x-e_ny_n|,
\end{multline}
for every $x\in B_{\delta\rho_0}(\delta\rho_0e_n)$,
$\delta$ small enough and depending on $E$ only.
By inequality \eqref{1N-13}, \eqref{star2}, \eqref{1-41},
\eqref{2N-11} and \eqref{1N-12} we obtain
\begin{multline}
\label{2N-13}
|\Gamma_+(x,t;y,0)-\Gamma_+(\Phi^{(t)}(x),t;y,0)|\\[2mm]
\leq C\left(\frac{(|x-e_ny_n|^2+t)}{\rho_0^2}\right)^{1/2}
\frac{\mathrm e^{-\frac{|x-e_ny_n|^2}{Ct}}}{t^{n/2}},
\end{multline}
for every $x\in B_{\delta\rho_0}(\delta\rho_0e_n)$, $y_n\in(-\delta\rho_0,0)$,
$t\in(0,(r_1/8)^2)$,
where $C$, $C\geq1$, $\delta$, $0<\delta<1$, depend on $E$ only.

We finally estimate
$|\nabla_x\Gamma(x,t;y,0)-\nabla_x\Gamma_+(x,t;y,0)|$
for $x\in B_{\delta\rho_0}(\delta\rho_0e_n)$, $y_n\in(-\delta\rho_0,0)$.
By \eqref{1N-3}, \eqref{1N-11}, \eqref{1-c5}, \eqref{2-c5} and Proposition \ref{5-16Tpr}
we have
\begin{multline}
\label{1N-14}
|\nabla_x\Gamma(x,t;y,0)-\nabla_x\Gamma_+(x,t;y,0)|\\
\leq \frac{C}{\rho_0}
\frac{\mathrm e^{-\frac{|x-e_ny_n|}{Ct}}}{t^{n/2}}
\left(\frac{[|x-e_ny_n|^2+t]^{1/2}}{\rho_0}\right)^{-1+\frac{\beta}{\beta+1}}\\
+C\left|\nabla_x\Gamma_+(x,t;y,0)-\nabla_x\left(\Gamma_+(\Phi^{(t)}(x),t;y,0)\right)\right|,
\end{multline}
for every $x\in B_{\delta\rho_0}(\delta\rho_0e_n)$, $y_n\in(-\delta\rho_0,0)$,
$t\in(0,(r_1/8)^2)$,
where $C$, $C\geq1$, depends on $E$ only.
Let us consider now the last term at the right hand side of \eqref{1N-14}.
We have
\begin{equation}
\label{2N-14}
\left|\nabla_x\Gamma_+(x,t;y,0)-\nabla_x\left(\Gamma_+(\Phi^{(t)}(x),t;y,0)\right)\right|
\leq H_1(x,t;y)+H_2(x,t;y),
\end{equation}
where
$$H_1(x,t;y)=C
|\nabla_x\Gamma_+(x,t;y,0)|\,|I-D_x\Phi^{(t)}(x)|$$
and
$$H_2(x,t;y)=C
|D_x\Phi^{(t)}(x)|\,|\left(\nabla_x\Gamma_+\right)(\Phi^{(t)}(x),t;y,0)-\nabla_x\Gamma_+(x,t;y,0)|,$$
where $C$ depends on $E$ only.
By Proposition \ref{5-16Tpr}, \eqref{e} and \eqref{N1-10} we have
\begin{equation}
\label{1N-15}
H_1(x,t,y)\leq C
|x-e_ny_n|\,
\frac{\mathrm e^{-\frac{|x-e_ny_n|^2}{Ct}}}{t^{\frac{n+1}{2}}}
\leq C'\frac{\mathrm e^{-\frac{|x-e_ny_n|^2}{2Ct}}}{t^{n/2}},
\end{equation}
for every $x\in B_{\delta\rho_0}(\delta\rho_0e_n)$, $y_n\in(-\delta\rho_0,0)$,
$t\in(0,(r_1/8)^2)$, where $C$, $C'$, depend on $E$ only.
To estimate from above the function $H_2$ we apply Theorem \ref{L}.
Let $\delta_1$ be defined as in Proposition \ref{fax} and let us denote
$$\rho=\frac{\delta_1}{4}[|x-e_ny_n|^2+t]^{1/2}.$$
By \eqref{d} and \eqref{N1-10} we have that there exists $\delta\in(0,1)$,
depending on $E$ only such that
\begin{equation}
\label{2N-15}
|\Phi^{(t)}(x)-x|\leq\frac{1}{2}\rho,
\end{equation}
for every $x\in B_{\delta\rho_0}(\delta\rho_0e_n)$.
Now Theorem \ref{L} and Proposition \ref{fax} yield
\begin{equation}
\label{1N-16}
[\nabla_x\Gamma_+(\cdot,t;y,0)]_{\beta;B'_\rho(x')\times(x_n^0,x_n+\rho)}
\leq C\rho^{-(2+\beta)}\frac{e^{-\frac{|x-e_ny_n|^2}{Ct}}}{t^{\frac{n-1}{2}}},
\end{equation}
where $x_n^0=\max\{0,x_n^0-\rho\}$ and $C$ depends on $E$ only,
$x\in B_{\delta\rho_0}(\delta\rho_0e_n)$, $y_n\in(-\delta\rho_0,0)$, $t\in(0,(\delta\rho_0)^2)$.
By \eqref{d}, \eqref{1N-16} and \eqref{2N-15} we obtain
\begin{equation}
\label{2N-16}
H_2(x,t;y)\leq\frac{C}{\rho_0}|x|^2\rho^{-(2+\beta)}\frac{e^{-\frac{|x-e_ny_n|^2}{Ct}}}{t^{\frac{n-1}{2}}},
\end{equation}
for every
$x\in B_{\delta\rho_0}(\delta\rho_0e_n)$, $y_n\in(-\delta\rho_0,0)$, $t\in(0,(\delta\rho_0)^2)$,
where $C$, $C\geq1$, depends on $E$ only.
Finally, \eqref{N1-10} and \eqref{2N-16} yield
$$H_2(x,t;y)\leq
\frac{C}{\rho_0}\frac{\mathrm e^{-\frac{|x-e_ny_n|^2}{Ct}}}{t^{n/2}}
\left(\frac{[|x-e_ny_n|^2+t]^{1/2}}{\rho_0}\right)^{-1+\beta}.$$
The last inequality, \eqref{1N-15}, \eqref{1N-14}, \eqref{2N-14} give
\begin{multline*}
\left|\nabla_x\Gamma(x,t;y,0)-\nabla_x\Gamma_+(x,t;y,0)\right|\\
\leq\frac{C}{\rho_0}\frac{\mathrm e^{-\frac{|x-e_ny_n|}{Ct}}}{t^{n/2}}
\left(\frac{[|x-e_ny_n|^2+t]^{1/2}}{\rho_0}\right)^{-1+\frac{\beta}{\beta+1}},
\end{multline*}
for every
$x\in B_{\delta\rho_0}(\delta\rho_0e_n)$, $y_n\in(-\delta\rho_0,0)$, $t\in(0,(\delta\rho_0)^2)$,
where $C$ depends on $E$ only.
\cvd

\subsection{Proof of Proposition \ref{prefb}}
\label{subsec-pro-prefb}
\noindent{\bf Proof of Proposition \ref{pr1}}
First of all, let us observe that
\begin{equation}
\label{(1)-5}
\Gamma_0(\xi,\tau;-\lambda_3 e_n,0)=
h^n\Gamma_0(h\xi,h^2\tau;-\lambda_3he_n,0)
\end{equation}
and
\begin{equation}
\label{(2)-5}
\Gamma_+^*(\xi,\tau;-\lambda_1 e_n,\lambda_2)=
h^n\Gamma_+^*(h\xi,h^2\tau;-\lambda_1he_n,\lambda_2h^2).
\end{equation}
Indeed  \eqref{(1)-5} is a trivial consequence of the definition of $\Gamma_0$.
Concerning \eqref{(2)-5}, it can be proved as follows.
Denote by $y^{(h)}=-\lambda_1he_n$, $s^{(h)}=\lambda_2h^2$ we have
\begin{equation}
\label{(3)-5}
\int_{\mathbb R^{n+1}}\left(\Gamma_+^*\partial_t\varphi
+(1+(k-1)\chi_+(x))\nabla_x\Gamma_+^*\cdot\nabla_x\varphi\right)dxdt
=\varphi(y^{(h)},s^{(h)}),
\end{equation}
for every $\varphi\in C^\infty_0(\mathbb R^{n+1})$.
In \eqref{(3)-5} $\Gamma_+^*$ denotes the function
$\Gamma_+^*(x,t;y^{(h)},s^{(h)})$.
Now for an arbitrary function $\psi\in C^\infty_0(\mathbb R^{n+1})$
put $\varphi(x,t):=\psi\left(\frac{x}{h},\frac{t}{h^2}\right)$ in \eqref{(3)-5}.
In the obtained integral we perform the change of variables $x=h\xi$, $t=h^2\tau$.
Thus, taking into account that $\chi_+(h\xi)=\chi_+(\xi)$, $h>0$, and denoting by
$W_h(\xi,\tau)$ the right-hand side of \eqref{(2)-5} we have
$$\int_{\mathbb R^{n+1}}\left(W_h\partial_\tau\psi+
(1+(k-1)\chi_+(\xi))\nabla_\xi W_h\cdot\nabla_\xi\psi\right)
d\xi d\tau=\psi(-\lambda_1 e_n,\lambda_2),$$
for every $\psi\in C_0^\infty(\mathbb R^{n+1})$.
Therefore
\begin{equation}
\label{(1)-6}
\partial_\tau W_h+\mathrm{div}((1+(k-1)\chi_+(\xi))\nabla_\xi W_h)
=-\delta(\xi+\lambda_1 e_n,\tau-\lambda_2)
\end{equation}
and, by the definition of $W_h$,
\begin{equation}
\label{(2)-6}
W_h(\cdot,\tau)=0\qquad\textrm{for every }\tau>\lambda_2.
\end{equation}
Finally, by the uniqueness for Cauchy problem \cite{ar},
by \eqref{(1)-6} and \eqref{(2)-6} we obtain \eqref{(2)-5}.

Now, performing the change of variable $x=h\xi$, $t=h^2\tau$
in the integral at the left-hand side of \eqref{(1)-4}, we get,
by \eqref{(1)-5} and \eqref{(2)-5},
\begin{equation}
\label{(1)-7}
I^{(h)}=h^{-n}I^{(1)}.
\end{equation}
Now, recall that
\begin{equation}
\label{(2)-7}
\mathcal F_{\zeta'}\left(\Gamma_0(\cdot,x_n,t;-\lambda_3 e_n,0)\right)
=\frac{\mathrm e^{-|\zeta'|^2t}}{\sqrt{4\pi t}}
\mathrm e^{-\frac{(x_n+\lambda_3)^2}{4t}}.
\end{equation}
In the case $k>1$, by Parseval formula, \eqref{(4)-2} and \eqref{(2)-7} we have
$$I^{(1)}=\left|\int_Y
M_{\lambda_2,\lambda_3}(\zeta',x_n,t;\rho)
F(\zeta',-\lambda_1;\rho)d\zeta'dx_ndtd\rho\right|,$$
where $Y=\mathbb R^{n-1}\times(0,+\infty)\times(0,\lambda_2)\times(0,1)$
and
\begin{multline*}
M_{\lambda_2,\lambda_3}(\zeta',x_n,t;\rho)=
\frac{|\zeta'|^2}{(2\pi)^{n-1}}
\frac{\mathrm e^{-t|\zeta'|^2}}{\sqrt{4\pi t}}
\mathrm e^{-\frac{(x_n+\lambda_3)^2}{4t}}\\
\times\left(|\zeta'|+\sqrt{\frac{k-1}{k}}\frac{(x_n+\lambda_3)}{2t}\sqrt\rho\right)
E(\zeta',x_n,\lambda_2-t;\rho).
\end{multline*}
For fixed $\lambda_2>0$, $\lambda_3>0$, taking into account
\eqref{(2)-2} and \eqref{(3)-2} we have
$$\lim\limits_{\lambda_1\to0^+}I^{(1)}=
\left|\int_Y M_{\lambda_2,\lambda_3}(\zeta',x_n,t;\rho)
\mathrm{Im}(A_1(\rho))d\zeta' dx_n dt d\rho\right|.$$
Thus
$$\lim\limits_{\lambda_1\to0^+}I^{(1)}>0$$
and by \eqref{(1)-7} the thesis follows.

Concerning the case $0<k<1$, we only give a sketch of the proof,
indeed such a case can be treated similarly to the case $k>1$.
In the case $0<k<1$, in order to have a suitable formula for $I^{(1)}$,
first we evaluate the integrals
\begin{eqnarray*}
&&\mathrm{Im}\left(\int_0^{+\infty}\mathrm e^{-i\mu x_n} K(x_n+\lambda_3,t)
dx_n\right),\\[2mm]
&&\mathrm{Im}\left(\int_0^{+\infty}\mathrm e^{-i\mu x_n}
\frac{\partial}{\partial x_n}K(x_n+\lambda_3,t) dx_n\right),
\end{eqnarray*}
where
$$\mu=\sqrt{\frac{1-k}{k}}\sqrt{1-\rho}|\zeta'|\qquad
\textrm{and}
\qquad K(x_n,t)=\frac{\mathrm e^{-\frac{x_n^2}{4t}}}{\sqrt{4\pi t}}.$$
In order to carry out such an evaluation we may use formula 3.322 of \cite{Gr}.
Then we choose $\lambda_2=\lambda_1^2$, we perform the change of variable
$\zeta'=\frac{\xi'}{\lambda_1}$, $t=\lambda_1^2\eta$ in the integral $I^{(1)}$
and we get
$$\lim\limits_{\lambda_1\to0^+}\left(\lambda_1^n
\lim\limits_{\lambda_3\to0^+} I^{(1)}\right)>0$$
and the thesis follows.
\cvd

\noindent{\bf Proof of Proposition \ref{prefb}.}
By the triangle inequality we get
\begin{eqnarray}
\label{(1)-13}
|\mathcal U(y_1,t_1;\overline y,\overline t)|
&=&|S_1(y_1,t_1;\overline y,\overline t)-S_2(y_1,t_1;\overline y,\overline t)|\\
&\geq&|S_1(y_1,t_1;\overline y,\overline t)|-|S_2(y_1,t_1;\overline y,\overline t)|.\nonumber
\end{eqnarray}
Let us first estimate from below $|S_1(y_1,t_1;\overline y,\overline t)|$.
Recall that $\Gamma_+^*(x,t;\overline y,\overline t)$ is the fundamental solution
of the adjoint operator of
$\mathcal L_+=\partial_t-\mathrm{div}((1+(k-1)\chi_+)\nabla)$.
Denote by
$$Q_{\rho/2}=B_{\rho/2}(0)\times(t_1,\overline t),
\quad Q^+_{\rho/2}=B^+_{\rho/2}(0)\times(t_1,\overline t),
\quad Q^-_{\rho/2}=B^-_{\rho/2}(0)\times(t_1,\overline t).$$
By the triangle inequality we have
\begin{equation}
\label{(1)-14}
|S_1(y_1,t_1;\overline y,\overline t)|\geq I_1-R_1-R_2,
\end{equation}
where
\begin{eqnarray}
\label{(2)-14}
&&I_1=\left|\int_{Q_{\rho/2}^+}\nabla_x\Gamma_+^*(x,t;\overline y,\overline t)
\cdot \nabla_x\Gamma_0(x,t;y_1,t_1)dxdt\right|,\\[2mm]
\label{(3)-14}
&&R_1=\int_{D_1(t_1,\overline t)\setminus Q_{\rho/2}}
|\nabla_x\Gamma_2(x,t;y_1,t_1)|\,|\nabla_x\Gamma_1^*(x,t;\overline y,\overline t)|dxdt\\[2mm]
&&\qquad +\int_{Q_{\rho/2}^+\setminus (D_1(t_1,\overline t)\cap Q^+_{\rho/2})}
|\nabla_x\Gamma_+^*(x,t;\overline y,\overline t)|\,|\nabla_x\Gamma_0(x,t;y_1,t_1)|dxdt,\nonumber
\end{eqnarray}
\begin{eqnarray}
\label{(4)-14}
&& R_2=\int_{D_1(t_1,\overline t)\cap Q_{\rho/2}}
|\nabla_x\Gamma_+^*(x,t;\overline y,\overline t)-\nabla_x\Gamma_1^*(x,t;\overline y,\overline t)|\\[2mm]
&&\qquad\qquad\times|\nabla\Gamma_0(x,t;y_1,t_1)|dxdt\nonumber\\[2mm]
&&\qquad +\int_{D_1(t_1,\overline t)\cap Q_{\rho/2}}
|\nabla_x\Gamma_0(x,t;y_1,t_1)-\nabla_x\Gamma_2(x,t;y_1,t_1)|\nonumber\\[2mm]
&&\qquad\qquad\times |\nabla_x\Gamma_1^*(x,t;\overline y,\overline t)|dxdt\nonumber,
\end{eqnarray}
where $D_1(t_1,\overline t)=\cup_{t\in(t_1,\overline t)}D_1(t)\times\{t\}$.

Now we estimate from below the term $I_1$.
First we notice that if $0<\delta\leq\frac{1}{4\sqrt 2}$ then
\begin{equation}
\label{(1)-15}
|x-\overline y|^2\geq\frac{1}{32}(|x|^2+\rho^2),\qquad
|x-y_1|^2\geq\frac{1}{32}(|x|^2+\rho^2),
\end{equation}
for every $x\in\mathbb R^n_+\setminus B_{\rho/2}^+$.
Also, we have trivially
\begin{multline}
\label{(2)-15}
I_1\geq
\left|\int_{\mathbb R^n_+\times(t_1,\overline t)}
\nabla_x\Gamma_+^*(x,t;\overline y,\overline t)
\cdot \nabla_x\Gamma_0(x,t;y_1,t_1)dxdt\right|\\
-\int_{(\mathbb R^n_+\times(t_1,\overline t))\setminus Q_{\rho/2}^+}
|\nabla_x\Gamma_+^*(x,t;\overline y,\overline t)|
\, |\nabla_x\Gamma_0(x,t;y_1,t_1)|dxdt.
\end{multline}
We now use Proposition \ref{5-16Tpr} and \eqref{(1)-15}
to estimate from above the second integral of the right-hand side of \eqref{(2)-15}.
We have
\begin{multline}
\label{(3)-15}
\int_{(\mathbb R^n_+\times(t_1,\overline t))\setminus Q_{\rho/2}^+}
|\nabla_x\Gamma_+^*(x,t;\overline y,\overline t)|
\, |\nabla_x\Gamma_0(x,t;y_1,t_1)|dxdt\\
\leq
C_0 \int_{\mathbb R^n_+\times(t_1,\overline t)}
\mathrm e^{-\frac{\rho^2}{C_0(t-t_1)}-\frac{\rho^2}{C_0(\overline t-t)}}
\frac{\mathrm e^{-\frac{|x|^2}{C_0(t-t_1)}-\frac{|x|^2}{C_0(\overline t-t)}}}
{(\overline t-t)^{\frac{n+1}{2}}(t-t_1)^{\frac{n+1}{2}}}
:=\tilde R,
\end{multline}
where $C_0$, $C_0\geq1$, depends on $k$ only.
Now performing the change of variables
$$z=\left(\frac{\overline t-t_1}{(t-t_1)(\overline t-t)}\right)^{1/2}x,$$
we have
\begin{multline*}
\tilde R\leq C_0\int_{\mathbb R^n\times(t_1,\overline t)}
\frac{\mathrm e^{-\frac{2\rho^2}{C_0(\overline t-t_1)}}}{(\overline t-t_1)^{n/2}}
\frac{\mathrm e^{-\frac{|z|^2}{C_0}}}{\sqrt{(\overline t-t)(t-t_1)}}dzdt\\
\leq\frac{C_3C_0}{\rho^n}\left(\int_{\mathbb R^n}\mathrm e^{-\frac{|z|^2}{C_0}}dz\right)
\left(\int_0^1\frac{d\lambda}{\sqrt{\lambda(1-\lambda)}}\right),
\end{multline*}
where $C_3=\max\limits_{s\in(0,+\infty)}\{s^{n/2}\mathrm e^{-2C_0s}\}$.
By the inequality obtained above, by \eqref{(2)-15}
and by Proposition \ref{pr1} we have
\begin{equation}
\label{(1)-16}
I_1\geq\frac{1}{C h^n}-\frac{C}{\rho^n},
\end{equation}
where $C$, $C\geq1$, depends on $k$ only.

In order to complete the proof we have to estimate from above the terms
$R_1$ and $R_2$ defined in \eqref{(3)-14}, \eqref{(4)-14}.
Denote by $R_{11}$ and $R_{12}$ the first and the second integral at the right-hand
side of \eqref{(3)-14} respectively. $R_{11}$ can be estimate in the same way of the integral
at the left-hand side of \eqref{(3)-15} and we have
\begin{equation}
\label{(1)-17}
R_{11}\leq\frac{C}{\rho^n},
\end{equation}
where $C$ depends on $k$ only.
Concerning $R_{12}$, by \eqref{q1} and Proposition \ref{5-16Tpr} we have
\begin{equation}
\label{(2)-17}
R_{12}\leq C_4\int_{t_1}^{\overline t}dt\int_{\mathbb R^{n-1}}dx'
\int_{-\psi(x',t)}^{\psi(x',t)}
\frac{\mathrm e^{-\frac{|x-y_1|^2}{C_4(t-t_1)}}}{(t-t_1)^{\frac{n+1}{2}}}
\frac{\mathrm e^{-\frac{|x-\overline y|^2}{C_4(\overline t-t)}}}{(\overline t-t)^{\frac{n+1}{2}}}
dx_n,
\end{equation}
where
$$\psi(x',t)=\frac{3}{2}\frac{E}{\rho_0}
\left(|x'|^2+|t-\overline t|\right)$$
and $C_4$, $C_4\geq1$, depends on $k$ only.
Now we perform, in the integral at the right-hand side of \eqref{(2)-17}
the following change of variables
\begin{eqnarray*}
t=t_1+\tau(\overline t-t_1),\quad x'=(\tau(1-\tau))^{1/2}z',
\quad x_n=(\tau(1-\tau))^{1/2}\xi-\lambda_1 h.
\end{eqnarray*}
Thus, denoting by
\begin{eqnarray*}
&&\sigma(\tau)=\frac{1}{\sqrt{\tau(1-\tau)}},\\[2mm]
&&\phi_1(z',\tau)=\frac{3}{2}\frac{E}{\rho_0}
\left(|z'|^2\sqrt{\lambda_2}\sqrt{\tau(1-\tau)}
+\frac{\sqrt{\lambda_2}\sqrt{1-\tau}}{\sqrt\tau}\right),\\[2mm]
&&A(z',\xi,\tau)=|z'|^2+\tau\xi^2+(1-\tau)
\left(\xi+\frac{(\lambda_3-\lambda_1)h}{\sqrt{\lambda_2}\sqrt{\tau(1-\tau)}}\right)^2,\\[2mm]
&&\theta(h)=\int_0^1\int_{\mathbb R^{n-1}}
\int_{\frac{\sigma(\tau)}{\sqrt{\lambda_2}}-h\phi_1(z',\tau)}^{\frac{\sigma(\tau)}{\sqrt{\lambda_2}}+h\phi_1(z',\tau)}
\frac{\mathrm e^{-\frac{A(z',\xi,\tau)}{C_4}}}{\sqrt{\tau(1-\tau)}}d\xi dz' d\tau,
\end{eqnarray*}
we get
\begin{equation}
\label{(1)-18}
R_{12}\leq\frac{C}{h^n}\theta(h),
\end{equation}
where $C$ depends on $k$ only.
Now observing that
$$A(z',\xi,\tau)\geq|z'|^2+\left(\xi+\frac{(\lambda_3-\lambda_1)}{\sqrt{\lambda_2}\sqrt\tau}\sqrt{1-\tau}\right)^2$$
and applying the H\"older inequality, we obtain
\begin{equation}
\label{(1)-19}
\theta(h)\leq\tilde C_p\left(\frac{h}{\rho_0}\right)^{1-\frac{1}{p}},
\end{equation}
for every $p\in(1,+\infty)$, where $\tilde C_p$ depends on
$p$ and $E$ only.
By \eqref{(1)-17}, \eqref{(1)-19} and recalling that $R_1=R_{11}+R_{12}$ we obtain
\begin{equation}
\label{(2)-19}
R_1\leq\frac{C}{\rho^n}+\tilde C_p\left(\frac{h}{\rho_0}\right)^{1-\frac{1}{p}}
\frac{1}{h^n},
\end{equation}
for every $h$, $0<h\leq\delta\min\{\rho,\sqrt{\overline t}\}$
and every
$\delta$, $0<\delta\leq\min\left\{\frac{\lambda_3}{C_1},\frac{1}{C_2},\frac{1}{4\sqrt 2}\right\}$,
where $C_1$ and $C_2$ are defined in \eqref{1-14T} and \eqref{(1)-12} respectively
and $C$ depends on $k$ only.

In order to estimate $R_2$, denote by $R_{21}$ and $R_{22}$ the first and the second integral
at the right-hand side of \eqref{(4)-14} respectively.
By Theorem \ref{STIMA ASINTOTICA} we have that there exists a constant $C_5$, $C_5\geq1$,
depending on $E$ only such that if $0<\delta\leq\frac{1}{C_5}$ and
$(x,t)\in\mathbb K_\rho:=\{(x,t)\in B_{\rho/C_5}\times(t_1,\overline t)\,:\,
x_n>\frac{1}{C_5\rho_0}(|x'|^2+|t-\overline t|)\}$ then
\begin{equation}
\label{(1)-20}
|\nabla_x\Gamma_1^*(x,t;\overline y,\overline t)-\nabla_x\Gamma_+^*(x,t;\overline y,\overline t)|
\leq \frac{C}{\rho_0^\alpha}
\frac{\mathrm e^{-\frac{|x-\overline y|^2}{C(\overline t-t)}}}
{(\overline t-t)^{\frac{n}{2}+\frac{1}{2}-\frac{\alpha}{2}}},
\end{equation}
where $C$, $C\geq1$, depends on $E$ and $k$ only and $\alpha=\frac{\beta}{\beta+1}$,
$\beta$ being defined in Theorem \ref{L}.
We have
\begin{equation}
\label{(2)-20}
R_{21}=J'+J'',
\end{equation}
where
\begin{eqnarray*}
&&J'=
\!\!\!\!\!\!
\int\limits_{D_1(t_1,\overline t)\cap Q_{\rho/2}\cap\mathbb K_\rho}
\!\!\!\!\!\!
|\nabla_x\Gamma_+^*(x,t;\overline y,\overline t)-\nabla_x\Gamma_1^*(x,t;\overline y,\overline t)|
\,|\nabla_x\Gamma_0(x,t;y_1,t_1)|dxdt,\\[2mm]
&&J''=
\!\!\!\!\!\!
\int\limits_{D_1(t_1,\overline t)\cap Q_{\rho/2}\setminus\mathbb K_\rho}
\!\!\!\!\!\!
|\nabla_x\Gamma_+^*(x,t;\overline y,\overline t)-\nabla_x\Gamma_1^*(x,t;\overline y,\overline t)|
\,|\nabla_x\Gamma_0(x,t;y_1,t_1)|dxdt.
\end{eqnarray*}
By \eqref{(1)-20} and Lemma \ref{fr} we have
\begin{equation}
\label{(1)-21}
J'\leq C\frac{h^{\alpha-n}}{\rho_0^\alpha},
\end{equation}
where $C$ depends on $E$ and $k$ only.
To estimate from above $J''$ we can arrange the method used to estimate $R_{11}$
and $R_{12}$ and we obtain that there exists $C_6\geq C_5$,
$C_6$ depending on $E$ and $k$ only
such that
if $0<\delta\leq\frac{1}{C_6}$ then, for $p\in(1,+\infty)$,
\begin{equation}
\label{(2)-21}
J''\leq\frac{C}{\rho^n}+\tilde C_p\left(\frac{h}{p}\right)^{1-\frac{1}{p}}
\frac{1}{h^n},
\end{equation}
where $C$ depends on $E$ and $k$ only and
$\tilde C_p$ depends on $p$ and $E$ only.
By choosing $p=\frac{1}{1-\alpha}$, \eqref{(2)-20}, \eqref{(1)-21}
and \eqref{(2)-21} yield
\begin{equation}
\label{(3)-21}
R_{21}\leq\frac{C}{\rho^n}+\tilde C\frac{h^{\alpha-n}}{\rho_0^\alpha},
\end{equation}
where $C$ depends on $E$ and $k$ only and $\tilde C$ depends on $p$ and $E$ only.

Now we estimate $R_{22}$. Denote by
$$w(x,t):=\Gamma_0(x,t;y_1,t_1)-\Gamma_2(x,t;y_1,t_1)$$
and recall \eqref{3-16T}.
We have that $w$ solves the heat equation in $B_{\rho/2}(y_1)\times(t_1,\overline t)$
and, since $w(x,t)=0$ for $(x,t)\in\mathbb R^n\times(-\infty,t_1]$
we can say that $w$ solves the heat equation in
$B_{\rho/2}(y_1)\times(\overline t-\rho^2,\overline t)$.
On $\partial B_{\rho/2}(y_1)\times(\overline t-\rho^2,\overline t)$
we have
$$|w(x,t)|\leq\frac{C}{(t-t_1)^{n/2}}
\mathrm e^{-\frac{\rho^2}{C(t-t_1)}}\chi_{[t_1,+\infty)}
\leq\frac{C'}{\rho^n}\chi_{[t_1,+\infty)},$$
where $C,C'$ depend on $k$ only.
Therefore by maximum principle and by standard regularity estimates we get
\begin{equation}
\label{(1)-22}
|\nabla_x w(x,t)|\leq\frac{C}{\rho^{n+1}}
\qquad\textrm{in }B_{\rho/4}(y_1)\times\left(\overline t-\frac{\rho^2}{4},\overline t\right).
\end{equation}
It is possible to have a similar estimate for $w$ in
$(\mathbb R^{n}\setminus B_{\rho/4}(y_1))\times(\overline t-\rho^2,\overline t)$,
namely by Proposition \ref{5-16Tpr} we have
\begin{equation}
\label{(2)-22}
|\nabla_x w(x,t)|\leq
\frac{C\mathrm e^{-\frac{\rho^2}{C((t-t_1)}}}{(t-t_1)^{\frac{n+1}{2}}}
\leq\frac{C'}{\rho^{n+1}},
\qquad (x,t)\in\left(\mathbb R^n\setminus B_{\rho/4}(y_1)\right)\times(\overline t-\rho^2,\overline t).
\end{equation}
By \eqref{(1)-22} and \eqref{(2)-22} we have
\begin{equation}
\label{N}
R_{22}\leq\frac{C}{\rho^{n+1}}\int_{t_1}^{\overline t}\int_{\mathbb R^n}
\frac{\mathrm e^{-\frac{|x-\overline y|^2}{\overline t-t}}}{(\overline t-t)^{\frac{n+1}{2}}}
dxdt\leq\frac{C'}{\rho^n},
\end{equation}
where $C,C'$ depend on $k$ only.

The estimate from above of $|S_2(y_1,t_1;\overline y,\overline t)|$ can be carried out
in a similar way of that used to estimate the integral in formula \eqref{(3)-15}.
Thus taking into account \eqref{3-16T}, we get there exists $C_7\geq C_6$ such that
if $0<\delta\leq\frac{1}{C_7}$ then
$$|S_2(y_1,t_1;\overline y,\overline t)|\leq\frac{C}{\rho^n},$$
where $C$ depends on $k$ only.
This inequality and \eqref{N}, \eqref{(3)-21}, \eqref{(2)-19}, \eqref{(1)-16},
\eqref{(1)-14}, \eqref{(1)-13} give \eqref{efb}.
\cvd

\end{document}